\documentclass{amsart}
\usepackage{graphicx}
\usepackage{amssymb}
\usepackage{amsmath, amsthm}
\usepackage{verbatim}
\usepackage[all]{xy}

\newcommand{\R}{\mathbb{R}}
\newcommand{\N}{\mathbb{N}}
\newcommand{\Z}{\mathbb{Z}}
\newcommand{\Heis}{\mathcal{H}}
\newcommand{\Pres}{\mathcal{P}}
\newcommand{\Qres}{\mathcal{Q}}
\newcommand{\hyp}{\mathbb{H}}
\newcommand{\VolD}{VolD}
\newcommand{\la}{\langle} 
\newcommand{\ra}{\rangle} 

\newcommand{\dist}{\textrm{VolD}}
\newcommand{\GL}{\textrm{GL}}

\theoremstyle{theorem}
\newtheorem{theorem}{Theorem}[subsection]

\newtheorem{cor}[theorem]{Corollary}

\newtheorem{lemma}[theorem]{Lemma}
\newtheorem{prop}[theorem]{Proposition}

\newtheorem{case}{Case}

\theoremstyle{remark}
\newtheorem{remark}{Remark}
\newtheorem{definition}{Definition}

\newtheorem{conjecture}{Conjecture}

\begin{document}

\title{Volume distortion in groups}
\author{Hanna Bennett}
\date{\today}

\begin{abstract}
Given a space $Y$ in $X$, a cycle in $Y$ may be filled with a chain in two ways: either by restricting the chain to $Y$ or by allowing it to be anywhere in $X$. When the pair $(G,H)$ acts on $(X, Y)$, we define the $k$-volume distortion function of $H$ in $G$ to measure the large-scale difference between the volumes of such fillings. We show that these functions are quasi-isometry invariants, and thus independent of the choice of spaces, and provide several bounds in terms of other group properties, such as Dehn functions. We also compute the volume distortion in a number of examples, including characterizing the $k$-volume distortion of $\Z^k$ in $\Z^k \rtimes_M \Z$, where $M$ is a diagonalizable matrix. We use this to prove a conjecture of Gersten.

\end{abstract}

\maketitle

\section{Introduction}

\subsection{Overview}

Consider a geodesic metric space $X$ with geodesic subspace $Y$. Given a pair of points in $Y$, there are two ways to measure the distance between them: we can consider the minimum of the lengths of paths between them that lie entirely in $Y$,  or we can allow paths to lie anywhere in $X$. Depending on how $Y$ is embedded in $X$, the latter distance may be much shorter. This idea can be generalized to higher dimensions: given a $(k-1)$-cycle $z$, we call the smallest volume of a $k$-chain whose boundary is $z$ the {\em filling volume} of $z$ and denote this $FV^k(z)$. When $z$ lies in $Y$, there are two possiblities: we may fill $z$ with a chain that lies anywhere in $X$, giving us $FV^k_X(z)$, or we might require that the chain be restricted to the subspace $Y$, which gives us $FV^k_Y(z)$. How do these two volumes compare? The {\em volume distortion function} provides a measurement of the difference in a large-scale sense.

We are particularly interested in the case in which $X$ and $Y$ are spaces on which a group and subgroup act cocompactly and properly discontinuously by isometries. We can always find such spaces by constructing a $K(G,1)$ CW complex that contains a $K(H,1)$ complex and then considering their universal covers. Counting $k$-cells gives us a combinatorial definition of $k$-volume in these spaces. If the Eilenberg-MacLane spaces have a finite $k$-skeleton, that is, $H$ and $G$ are $F_{k}$, then we may speak of \textit{subgroup volume distortion}, by which we mean volume distortion in the spaces on which the groups act.

\begin{definition}
Let $H$ be a subgroup of $G$, both $F_{k}$ groups, and let $X$ be the universal cover of an Eilenberg-MacLane space of $G$ and  $Y \subset X$ the universal cover of an Eilenberg-MacLane spaces for $H$.  The \textit{$k$-volume distortion function} function of $H$ in $G$ is a function $\VolD^{k}_{(G, H)}: \N \to \N$ given by

$$\VolD^{k}_{G, H}{(n)} = \max \{ FV^{k}_Y(z) \mid z \textup{ is a } (k-1)-\textup{cycle in } Y \textup{ and } FV_X(z) \leq n\}$$
\end{definition}

Notice that if the filling volume is the same in the subspace as the ambient space, we get a linear volume distortion function. Thus we say that a subgroup $H$ is {\em$k$-volume undistorted} in $G$ if the volume distortion function is linear.

While length distortion is well-understood and area distortion has been studied to some extent (see \cite{gersten96}), higher-dimensional volume distortion is new.

We prove a number of foundational facts in \S \ref{theorems}: up to linear terms, the distortion functions of two pairs of quasi-isometric CW-complexes are equivalent (Theorem \ref{qi}), and thus that volume distortion is independent of the choice of spaces. In this section we also provide bounds in terms of $k^{th}$-order Dehn functions and discuss the computability of volume distortion functions. We then compute a number of examples in \S \ref{exs}.

In \cite{gersten96}, Gersten proves that the copy of $\Z^2$ is always area-undistorted in $\Z^2 \rtimes_M \Z$ (note that here $M \in \GL(2, \Z)$). He gives the following conjecture. 

\begin{conjecture}[Gersten, \cite{gersten96}, p. 19]\label{conjgersten}
The group $\Z^k$, $k \geq 3$, is area undistorted in $\Z^k \rtimes_M \Z$ if and only if $M$ is of finite order in $GL(k,\Z)$. 
\end{conjecture}

In \S \ref{gerstenconj} we prove a generalization of this conjecture; we allow $M$ to be any $m$-by-$m$ integer-entry matrix, and consider the group
$$\Gamma_M = \la x_1, \cdots x_m, t  \mid [x_i, x_j] = 1, tx_it^{-1} = \phi(x_i) \text{ for } 1 \leq i, j \leq m \ra,$$ where $\phi$ is a homomorphism taking $x_i$ to $x_1^{a_1}x_2^{a_2} \cdots x_m^{a_m}$, where the $a_j$ form the $i$th column of $M$. When $\det M = 1$, we can write $\Gamma_M$ as the semidirect product $\Z^m \rtimes_M \Z$. We then prove the following theorem.

\begin{theorem}\label{areadistthm}
$\Z^m$ is area-undistorted in $\Gamma_M$ if and only if $M$ has finite order.
\end{theorem}

Conjecture \ref{conjgersten} is then the special case when $\det M = 1$. 

\begin{center}
\begin{figure}
\includegraphics[scale=.6]{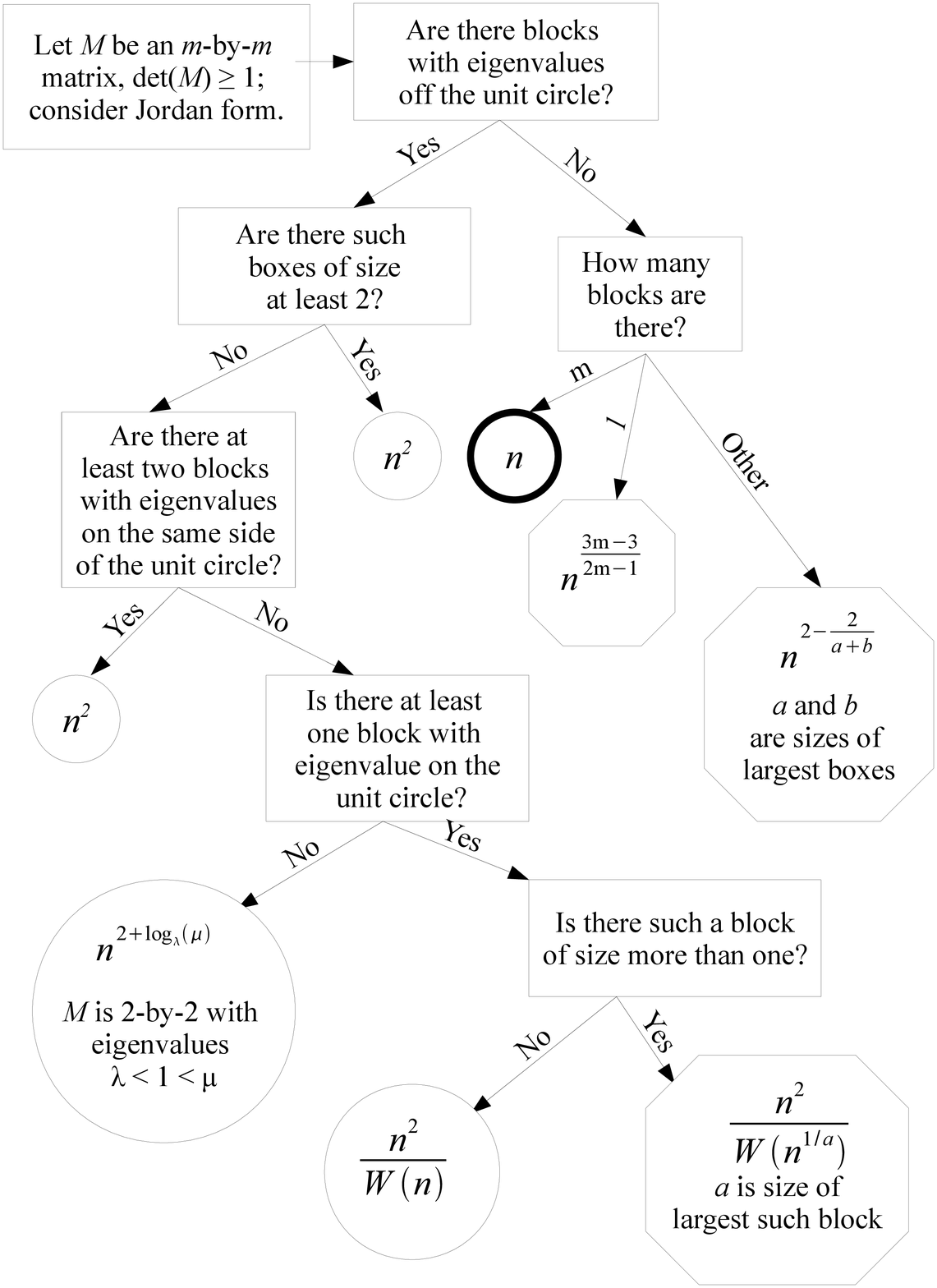}
\caption{Area distortion flow chart. Circles indicate sharp bounds while octogons are lower bounds only; $n^2$ is always an upper bound. The function $W$ is the Lambert $W$ function, that is, the inverse of $xe^x$.}\label{areaflow}
\end{figure}
\end{center}

Theorem \ref{areadistthm} is proved by identifying different cases and calculating a lower bound for the area distortion function in each case. This is illustrated in Figure \ref{areaflow}, which charts the possible cases and the resulting area distortion in each case.

We generalize this to look at higher volume distortion of $\Z^k$ in $\Gamma_M$. In the case of diagonalizable matrices, we can characterize the $k$-volume distortion completely.

\begin{theorem}\label{gammamdist}
Let $M$ be an integer-entry $k$-by-$k$ diagonalizable matrix with
 $det(M) = d \geq 1$,and let $\lambda_i$ denote the absolute value of the $i^{th}$ eigenvalue. Then the $k$-volume distortion of $\Z^k$ in $\Gamma_M$ depends only on the eigenvalues of $M$. If $M$ has at least two eigenvalues off the unit circle, the volume distortion is 
$$\dist^{(k)}(n) \asymp n^{1 + \log{d}/\log{\alpha}}, \textrm{ where } \alpha = \prod_{i = 1}^{k} \min \{\lambda_i, d\})/d.$$

If $M$ has exactly one eigenvalue off the unit circle, 
$$ \dist^{(k)}{n} \asymp (\frac{n^{k}}{W(n)})^{1/(k-1)}.$$

Otherwise, $\dist^{(k)}(n) \asymp n$.

\end{theorem}

Note that here we obtain a sharp bound. When $M$ is not diagonalizable, we provide a lower bound for volume distortion in \S \ref{nondiag}.

We can look more generally at groups of the form $G = H \rtimes_\phi \la t\ra$, where $H$ is any group and $\phi$ is an automorphism on $H$. There is a natural surjection to $\Z$, given by the second factor, that allows us a well-defined notion of {\em height} in the group and in a $\widetilde{K(G,1)}$ such that the height zero subspace is a $\widetilde{K(H,1)}$. Then we can think of $\phi$ as sending a $k$-cell at height $h$ to its image under $\phi$ at height $h$-1. This corresponds to conjugating by $t$ in the group presentation. 

The dynamical properties of $\phi$ can thus be used to find bounds on the volume distortion of $H$.  In \S \ref{cxitysection} we take the idea of complexity from \cite{gersten96} and alter it, so that the {\em complexity} of $\phi$, denoted $c_k(\phi)$, is the maximal $k$-volume of the image of a $k$-cell.

\begin{theorem}
\label{complexity}
Let $\phi$ be an automorphism on $H$, an $F_{k}$ group, and\\
 $m = \max \{c_k(\phi), c_{k}(\phi^{-1}\}$. Then $\VolD^{k}_{( H \rtimes \Z, H)}(n) \leq n \cdot m^n$. 
\end{theorem}

\begin{cor}\label{cxitycor}
When $\phi$ has complexity $m = 1$, then $G$ is $k$-undistorted in $G \rtimes_\phi \Z$.
\end{cor}

In particular, if a $K(G,1)$ has only one $k$-cell, then $\phi$ must send this $k$-cell to itself, because it induces an automorphism on the $k$-skeleton. Thus the $k$-complexity is one, so $G$ is $k$-undistorted in $G \rtimes_\phi \Z$.

\subsection{Acknowledgements}

I would like to thank my advisor, Benson Farb, for all his input and support; Noel Brady and Max Forester, for providing a key step to the proof of Theorem \ref{abc}; and everyone else whose conversations helped make this paper what it is, including Shmuel Weinberger, Robert Young, Nathan Broaddus, Irine Peng, and Tim Riley.

\section{Background}

Let $G$ be a group with presentation $\la S \mid R \ra$. We say a word in $F(S)$, the free group on $S$, is {\em null-homotopic} if represents the identity in $G$, that is, it can be written as a product of conjugates of relators. The {\em area} of a null-homotopic word is the minimal number of such relators necessary. The Dehn function for $G$, denoted by $\delta: \N \to \N$, is defined by
$$\delta(n) = \{ A(w) \mid l(w) \leq n\}$$
which provides an upper bound on the area of a word in terms of its length. While this function appears to depend on the presentation, we can create a relation of functions $f \preceq g$ when there exists some $C > 0$ such that 
\begin{equation}\label{fneq}
f(x) \leq Cg(Cx+C) + Cx + C
\end{equation}
and we say $f \asymp g$ if $f \preceq g$ and $g \preceq f$. Under this equivalence, the Dehn function is a quasi-isometry invariant, and so in particular independent of presentation.

This function can be used to answer the {\em word problem}, first asked by Dehn in \cite{dehn11}: given a word in $F(S)$, is there an algorithm for determining whether this word represents the identity? The answer is yes if and only if the Dehn function is computable. However, the algorithm provided by the Dehn function may not be very efficient. For example, if a group with exponential Dehn function can be embedded in a group with quadratic Dehn function, we can use the ambient group to more easily solve the word problem in the subgroup. In such a case, we may think of the embedding as being (area) distorted. 

\subsection{Definitions}
\subsubsection{Area distortion}

In \cite{gersten96}, Gersten defines a function, similar to the Dehn function, which measures this area distortion. Precisely, let G be a group with finite presentation $\Pres = \la S \mid R \ra$, and let $H$ be a subgroup with presentation $\Qres = \la S' \mid R' \ra$, where $\Qres$ is a subpresentation of $\Pres$, that is, $S' \subset S$ and $R' \subset R$. Then the {\em area distortion function of $H$ in $G$}, $Ad:\N \to \N$, is given by 
\begin{displaymath}
Ad(n) = \max \{Area_H(w) \mid Area_G(w) \leq n, w \in N(R')\}.
\end{displaymath}

It is not \emph{a priori} clear that such a maximum must exist---perhaps we could find a sequence of words representing the identity in $H$ with area in $G$ bounded by $n$, but area in $H$ growing arbitrarily large. This, however, cannot happen, precisely because $G$ and $H$ are finitely presented.

\begin{prop}[Gersten] The area distortion function is well-defined.

\end{prop}

\begin{proof}
Let $m$ be the length of the longest relator in $R$, and let $w$ be a word with van Kampen diagram of area at most $n$. Separate the diagram into a collection of topological circles, each with area $n_1, n_2, ..., n_k$ (note that $k$ and the $n_i$ are all bounded above by $n$). In such a topological circle of area $n_i$, the length of the boundary cannot be more than $m \cdot n_i$. Thus there is a finite number of possible loops for the topological circle; for each, we can fill in $H$ with some area. Combining these for each $i$ gives an upper bound on the area in $H$ of the word $w$ which depends only on $n$. Thus $AD(n)$ is bounded above for each $n$, so the function is well-defined.
\end{proof}

Note the importance here of dividing $w$ into pieces that contribute to the area. While we cannot bound the length of $w$, bounding the lengths of these pieces will often suffice for our purposes. We will continue to use this approach to bounding volume distortion, and so it will benefit us to give a name to the boundary of the ``area-contributing'' pieces of $w$. Let $D$ be a van Kampen diagram for $w$. Define the {\em frontier} of $D$ by $FR(D) = \partial (D^\circ)$. (Note that generally frontier is used as a synonym for boundary; we are modifying the definition to a subset of the boundary that will play an important role in bounding volumes.)

We will use the same equivalence for distortion functions as used for Dehn functions, given in Equation~\ref{fneq}. If $Ad$ is linear, then we say that the area of $H$ is {\em undistorted} in $G$, as this means that there is essentially no advantage to filling in $G$ over restricting to $H$. While this function is closely related to the Dehn functions of both the group and subgroup, the distortion function cannot in general be written simply as some combination of the Dehn functions of $H$ and $G$.

 Another related concept is \textit{length distortion}, often called simply {\em subgroup distortion}, which compares the lengths of elements in the subgroup to the lengths in the ambient group. These two concepts are independent: groups may have distorted length but undistored area (e.g. Sol groups), or undistorted length and distorted area (e.g. examples constructed in \cite{bbms}). 

\subsubsection{Volume distortion}

Just as Dehn functions have been generalized to higher dimensions (as 
``higher order'' Dehn functions $\delta^{(k)}$), we would like to generalize to volume distortion. To do so, we will need to take a more geometric approach than the algebraic one used for area distortion. We will first define volume distortion on CW-complexes, and then define volume distortion for groups in terms of complexes on which they act.

Because we are using pairs of groups, we must take some care to specify that the action respects this pairing. We say the pair $(G,H)$ acts geometrically on $(X,Y)$ if $G$ acts geometrically on $X$ (that is, $G$ acts cocompactly and properly discontinuously by isometries), and we can restrict this to an action of $H$ on $Y$ (if $h \in H$ and  $y \in Y$, then $ h.y \in Y$), and this is also a geometric action by isometries. Note that in particular we may construct a $K(H,1)$ inside of a $K(G,1)$ so that $(G, H)$ will act geometrically on $(\widetilde{K(H,1)}, \widetilde{K(G,1)})$; often this is what we will be considering. 

In order for the distortion function to be well-defined, we need to put some conditions on the CW-complex $X$. The conditions needed are exactly those given in \S 3 of \cite{alonsowangpride99} for a $k$-Dehn complex. $X$ is {\em k-Dehn} if $X$ is $k$-connected, the $m$-order Dehn functions are well-defined for $m \leq k$, and there is a uniform bound, say $r$, on the number of faces on an $m$-cell, for $m \leq k+1$.

In such a space $X$, a cellular $k$-chain is denoted by $z = \sum \alpha_i \sigma_i$ where the $\alpha_i$ are integers and $\sigma_i$ are $k$-cells. The \emph{volume} of $z$ is $V^k(z) = \sum |\alpha_i|$. Given a $k$-cycle $z$, we define the \emph{filling volume} of $z$, $FV^{k+1} (z)$, to be minimal volume over all $k$-chains which extend $z$, that is, 
$$FV^{k+1}(z) = \min \{V^{k+1} (u) \mid \partial u = z\}.$$ Since $X$ is $k$-connected, every cycle $z$ is the boundary of some chain. Note that this is the definition given in chapter 10 of \cite{echlpt}, but what we call volume they call mass, and what we call filling volume, they call volume. 

Given a subcomplex $Y$ of $X$, the \emph{k-volume distortion function} $\dist^k_{(X,Y)} : \N \to \N$ is defined by 
\begin{equation}
\dist^k(n) = \max \{FV^{k}_Y(z) \mid FV^k_X(z) \leq n, \text{ where $z$ is a $(k$--$1)$-cycle in $Y$} \}
\end{equation}

 The uniform bound on the size of the boundary of an $m$-cell serves the same purpose as the finite presentation in the definition of $AD$, that is, it ensures that a maximum exists. In particular, note that while we have no bound on the volume of $z$, we do obtain $rn$ as a bound on the volume of the frontier of any filling of $z$.

Let $(G,H)$ act as a pair on $(X, Y)$ by isometries, where the action is cellular and discrete. As mentioned above, $Y = \widetilde{K(H, 1)}$ and $X = \widetilde{K(G,1)}$ with $Y$ in $X$, will satisfy this. Then $\dist^k_{(G, H)} (n) := \dist^k_{(X, Y)}(n)$. This definition only makes sense if the spaces are $k$-Dehn; this will happen when $H$ and $G$ are $F_k$, that is, their $K(\pi,1)$'s have finite $k$-skeleton.

Note that $AD$ and $\dist^2$ are actually different functions: in the former case, the function involves homotopy, while in the latter the function involves homology. $\dist^2$ is referred to as \emph{weak distortion} in \cite{gersten96}. Both homology and homotopy definitions exist for higher-dimensional Dehn functions; while they are often the same in examples in which they are easy to compute, they are not equivalent in general. 

\subsubsection{Riemannian manifolds}

Sometimes we will be able to determine the volume distortion by considering group actions on Riemannian manifolds. In a general sense, the function will work the same way; what changes is the manner in which we define the volume and filling volume.

Let $M$ be a connected Riemannian manifold with submanifold $N$, such that a pair of groups $(G, H)$ act on $(M, N)$ properly discontinuously by isometries. We will work with Lipschitz $K$-chains, that is, formal finite sums with coefficients in $\{+1, -1\}$ of maps $f_i:\triangle_k \to M$ where $f_i$ is $K$-lipschitz for some universally fixed $K$. We choose lipschitz maps so that the functions are differentiable almost everywhere, leading to a well-defined idea of volume, and so that under quasi-isometry the composition with a lipschitz map is a bounded distance away from a lipschitz map.

We find the volume of a $k$-chain in the following way: for each lipschitz map $f$, we can consider $D_xf$ at almost every point in the domain. This map sends an orthonormal basis in $T_x \triangle_k$ to a set of vectors in $T_{f(x)}M$. These vectors give a parallelopiped; call its volume V(x). This is the $k$-dimensional Jacobian of $f$ at $x$, and can be found by considering the matrix
\begin{displaymath} A = \left(\begin{array}{c} D_xf(e_1) \\   D_xf(e_2) \\  \vdots \\  D_xf(e_k) \end{array}\right) 
\end{displaymath}
and taking $V(x) = \sqrt{det[A \cdot g \cdot A^T]}$, where $g = (g_{i,j})$ is the Riemannian metric. 

Now we integrate over $\triangle_k$:
$$V(f) = \int_{\triangle_k} V(x) dx$$

More generally, the volume of a $k$-chain is the sum of the volumes of the component maps $f_i$.

With this new definition of volume, we may now define the filling volume and distortion function just as before: given a lipschitz $(k-1)$-cycle $z$,  
\begin{equation}
FV^{(k)}(z) = \inf \{V(u) \mid u \text{ is a lipschitz $k$-chain with } \partial u = z\}
\end{equation}

Before we can define a volume distortion function in this case, we need one last requirement. Choose some $c \in \R^+$. Just as the value of $K$ does not matter so long as our maps are $K$-lipschitz for some $K$, this choice of $c$ will not affect the distortion function up to the usual equivalence of functions, which allows us to discuss ``the'' distortion function without specifying $c$ or $K$. 

The $k$-volume distortion function (with respect to c) is a function $\dist_{(M,N)}^{(k)}:\N \to \N$ with 
$$\dist(n) = \sup\{FV^{(k)}_N(z)\; \mid \; \exists \; k \textrm{-chain } u \text{ with } V^{k}_M(u) \leq n, FR(u) \leq cn \}.$$ 

Note that, in the case of CW-complexes with a cocompact group action, there is a natural choice of $c$: the maximal boundary volume of a $k$-cell. Once again, the restriction on the size of the frontier of a filling gives this function an upper bound (in terms of the Dehn function of the subspace), so a supremum exists.

In the following section we will show that the volume distortion in this case is equivalent to the version obtained by taking a triangulation that is invariant under $(G,H)$.

\section{General theory}

\label{theorems}

\subsection{Equivalence of definitions}

The definition of volume distortion appears to depend on the choice of spaces; however, we will show that up the equivalence of functions given above, it is a quasi-isometry invariant, and therefore in particular a group invariant. 

A \emph{$(K, C)$-quasi-isometric embedding} is a map $f: X \to X'$ where $$\frac{1}{K}d(x,y) - C \leq d(f(x), f(y)) \leq K d(x, y) + C$$ for all $x, y \in X$. This function is a quasi-isometry given the additional requirement that for all $y \in X'$ there is some $x \in X$ with $d(y, f(x)) \leq C$. This is equivalent to saying that $f$ has a $(K, C)$-quasi-isometric inverse  $g:X' \to X$ with $d(x, g(f(x))) \leq C$ and $d(x, f(g(x))) \leq C$.  

As with group actions, we will be considering pairs of spaces $(X, Y)$ where $Y \subset X$. We will say the pairs of spaces $(X_1, Y_1)$ and $(X_2, Y_2)$ are qusi-isometric if there is a quasi-isometry $f: X_1 \to X_2$ with $f(Y_1) \subseteq Y_2$ and $f|_{Y_1}: Y_1 \to Y_2$ also a quasi-isometry.\\

\begin{theorem}\label{qi}Suppose that $(X_1, Y_1)$ and $(X_2, Y_2)$ are $k$-Dehn spaces which are quasi-isometric as pairs. Then the distortion functions $\dist^k_1$ of $(X_1, Y_1)$ and $\dist^k_2$ of $(X_2, Y_2)$ are equivalent.
\end{theorem}
\begin{proof}

Suppose we have $(K, C)$-quasi-isometries $f:(X_1, Y_1) \to (X_2, Y_2)$ and quasi-isometric inverse $g:(X_2, Y_2) \to (X_1, Y_1)$.

See Figure ~\ref{qiillo} for an illustration of this proof. We start, as in Figure \ref{qiillo}(a), with a $(k-1)$-cycle $z$ in $Y_1$, which is filled in $X_1$ with a $k$-chain $u$, where $V(u) = n$. We want to fill $z$ in $Y_1$ with volume linear in $n$ and $\dist^k_2$. 

\begin{figure}[htbc]
\begin{center}
\includegraphics[width=5in]{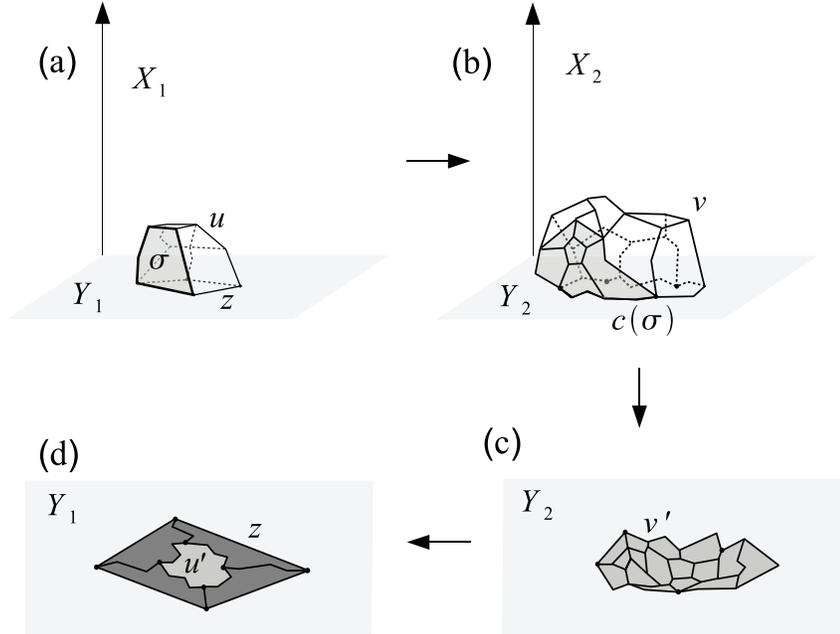}
\end{center}
\caption[Quasi-isometries preserve volume distortion]{Quasi-isometries preserve volume distortion}\label{qiillo}
\end{figure}

In Figure \ref{qiillo}(b) we  construct a $k$-chain $v$ in $X_2$ using $f(u^{(0)})$, that is, the image of the vertices in $u$. We want $v$ to have boundary in $Y_2$ and volume at most $A_kn$, where $A_k$ is independent of $u$. We do this by noting that for any $m$-cell $\sigma$ in $X_1$, we can construct an $m$-chain $c(\sigma)$ in $X_2$ with $m$-volume bounded by some constant, say $A_m$. Further, if the $m$-cell was originally in $Y_1$, then we can construct the new $m$-chain in $Y_2$ with the same bound $A$ on the volume. We do this construction by induction on $m$: when $m = 1$, we simply note that the distance between the image of the boundary points of a cell are at a distance at most $K+C$ apart, and choose a geodesic between these points; this gives a 1-chain with length at most $C+K$. 

Now suppose such a construction exists for dimension $(m$-$1)$ and let $\sigma$ be an $m$-cell. Carry out the construction on the boundary of $\sigma$. Since the spaces are $k$-Dehn, there is a universal bound $r$ on the volume of the boundary of any $m$-cell in $X_1$ when $m \leq k+1$. Thus we have constructed a $(k-1)$-cycle of $(k-1)$-volume at most $A_{k-1}r$. Again because the spaces are $k$-Dehn, we know that the $k^{th}$-order Dehn fuction of $X_2$ is well-defined, so the cycle can be filled in $X_2$ with volume at most $\delta_{X_2}^{(k)}(A_{k-1}r)$. If the original $k$-cell is in $Y_1$, then we can do the same thing, now filling with volume bounded by $\delta_{Y_2}^{(k)}(A_{k-1}r)$. The maximum of these two numbers is then the value of $A_k$. This constant depends on the spaces, quasi-isometries, and $k$, but is independent of $u$.

Using this procedure on each $k$-cell in $u$, we construct a $k$-chain $v$ in $X_2$ with boundary in $Y_2$ and $k$-volume at most $A_kn$. By definition, the cycle $\partial v$ can be filled in $Y_2$ with $k$-volume at most $\dist^k_2(A_kn)$. Call the chain with this filling $v'$; note that this chain lies entirely in $Y_2$. This is illustrated in Figure \ref{qiillo}(c). Just as above we can now construct a $k$-chain $u'$ in $Y_1$ which contains the 0-skeleton of $g(v')$ and has with volume at most $B_k\cdot \dist^k_2(A_kn)$, again with $B_k$ independent of $u$. The boundary of this chain will probably not be $z$, but each vertex in $z$ will correspond to a vertex in $u'$ which is a distance at most $C$ from the vertex in $z$; we can use this fact to construct a homotopy between the two cycles $z$ and $\partial u'$. Again we do this by building up by dimension on each cell of $z$. First we construct paths of length at most $C$ between a vertex in $z$ and the corresponding vertex in $u'$. Given two vertices in $z$ connected by an edge, their corresponding vertices are connected by a path of length at most $K(K+C)+C$. This gives us a loop of length at most $K(K+C) + 3C + 1$; this can be filled in with area at most $\delta_{Y_1}(K(K+C) + 3C + 1)$. See Figure \ref{qiillo}(d).

Continue this construction one dimension at a time: Given an $m$-cell in $z$, we can construct a homotopy between each of the boundary cells and a corresponding $m-1$-chain in $u'$ with some volume bounded by $D_{m-1}$, again, agreeing on their boundaries. This gives an $m$-cycle with volume at most $1 + rD_{m-1} + A_{m}B_m$, which can be filled with $m$-volume at most $D_k = \delta_{Y_1}^{(m)}(1 + rD_{m-1} + A_{m}B_m)$. Ultimately this gives us a homotopy between $z$ and $\partial u'$ with $k$-volume at most $D_kV(z)$. Since $V(z) \leq rn$, combining this with $u'$ gives us a filling of $u$ with volume bounded above by $D_krn + B_k\dist_2^k(A_kn)$. Thus $\dist_1^k(n) \leq D_krn + B_k\dist_2^k(A_kn)$. We can use the same process to reverse the roles of $\dist_1^k$ and $\dist_2^k$, so the two functions are equivalent.
\end{proof}

\begin{theorem}\label{defeq}
Let $M$ be a $(k-1)$-connected Riemannian manifold, with $(k-1)$-connected submanifold $N$. Let the pair $(G, H)$ be $F_k$ groups, that is, have $K(\pi, 1)$ with finite $k$-skeleton, and suppose they act properly discontinuously, cocompactly and by isometries as a pair on $(M, N)$. Let $\tau$ be a $G$-invariant triangulation of $M$, with $\tau_1$ an $H$-invariant subtrangulation of $\tau|_N$. Then the distortion function of $H$ in $G$, which by definition is the distortion function $\dist_\tau^k$ of the $k$-skeleton of $\tau_1$ in $\tau$, is equivalent to the geometric distortion function $\dist^k_{(M, N)}$.

\end{theorem}

\begin{proof}
The work done by Burillo and Taback in \cite{burillotaback02} to show the isoperimetric version of this theorem provides us with all the tools we need for this proof. In particular, they prove that, given $M$, $G$, and $\tau$ as above, the following holds.

\begin{lemma}[Pushing Lemma, Lemma 2.1 of \cite{burillotaback02}]\label{bt1} There exists a constant $C$, depending only on $M$ and $\tau$, with the following property: Let $T$ be a Lipschitz $(k-1)$-chain in M, such that $\partial T$ is included in $\tau^{(k-2)}$. Then there exists another Lipschitz $(k-1)$-chain $R$, with $\partial R = \partial T$, which is included in $\tau^{(k-1)}$, and a Lipschitz $k$-chain $S$, with $\partial S = T - R$, satisfying $V(R) \leq C V(T)$ and $V(S) \leq C V(T).$
\end{lemma}

Essentially, this is saying that $(k-1)$-chains in $M$ are very near $(k-1)$-chains in $\tau^{(k-1)}$ of comparable area.

Now let $z$ be a $(k-1)$-cycle in $N$, and let $u$ be a $k$-chain in $M$ with $\partial u = z$ and $V^k(u) \leq n$. Because $z$ has no boundary, we can apply Lemma~\ref{bt1} to get a $(k-1)$-cycle $z'$ in $\tau_1$, with homotopy between them given by the chain $S$ in $N$, where $V^k(S) \leq Cn$. Now $z'$ is a cycle in $\tau_1$; we can fill it with the chain $u - S$, so $FV(z') \leq (C+1)n$. Next we can apply Lemma~\ref{bt1} to $u - S$; this gives us a chain in $\tau$ with volume at most $C(C+1)n$ with boundary $z'$. By definition, we can fill $z'$ with some chain in $\tau_1$ with volume at most $\dist_\tau(C(C+1)n)$. Since $G$ acts cocompactly on $\tau$, there is some maximal volume, say $A$, of any $k$-cell in $\tau$; so we can now fill $z'$, say by $u'$, in $N$ with volume at most $A\dist_\tau(C(C+1)n)$. Then $u' - S$ fills $z$ with volume at most $A\dist_\tau(C(C+1)n) + Cn$.

On the other hand, suppose $z$ is a cycle in $\tau_1$, with a filling $u$ in $\tau$ which has volume at most $n$. Then we can fill $z$ in $M$ with volume at most $An$. This means we can fill $z$ in $N$ with a chain $u$ with volume at most $\dist_{(M,N)}(An)$. By applying \ref{bt1}, we can find a filling $u'$ of $z$ in $\tau_1$ with volume at most $C\dist_{(M, N)}(An)$.

\end{proof}

Note that the above proof is independent of the $c$ used to bound the volume of the boundary in the Riemannian manifold case; thus this constant does not affect the distortion function.

By combining Theorems \ref{qi} and \ref{defeq}, we obtain the following.

\begin{theorem}\label{riemeq}
Given pairs of spaces $(M_1, N_1)$ and $(M_2, N_2)$ which are quasi-isometric, and groups $(G_1, H_1)$ and $(G_2, H_2)$ where $(G_i, H_i)$ act cocompactly and properly discontinuously by isometries on $(M_i, N_i)$, the distortion function\\ $\dist^k_{(M_1,N_1)}$ is equivalent to the function $\dist^k_{(M_2, N_2)}$. 
\end{theorem}

\begin{proof}

Because the group actions are geometric, each pair $(G_i, H_i)$ is quasi-isometric to $(M_i, N_i)$. Since the $(M_i, N_i)$ are quasi-isometric to each other by assumption, the $(G_i, H_i)$ must be as well. Thus by Theorem~\ref{qi} their distortion functions are equivalent. By Theorem~\ref{riemeq}, the distortion functions of the pairs of Riemannian spaces are equivalent to those of the respective pair of groups. Thus the distortion functions of the spaces are equivalent as well.
\end{proof}

\subsection{Distortion and Dehn functions}
\label{dehn}

Because of the closely related definitions of Dehn functions and volume distortion functions, it can be tempting to believe that one can express the distortion function easily in terms of the Dehn functions of the group and subgroup. The reality is not so simple---for example, the subgroup may have a greater or smaller Dehn function than the ambient group. However, we can use Dehn functions to provide certain bounds for volume distortion functions.


\begin{theorem}
\label{dehnupperbound}
Let $H \subset G$ be $F_k$ groups, and let $\delta^{(k-1)}_H$ be the $(k-1)$-order Dehn function of $H$. Then $\dist_{(G, H)} \preceq \delta^{(k-1)}_H$.
\end{theorem}

\begin{proof}
Suppose $z$ is a $(k-1)$-cycle in $H$, with filling in $G$ of $k$-volume $n$. Because these groups are $F_k$, there is an upper bound on the number of boundary faces in a $k$-cell in $G$, say $r$. Then the $(k-1)$-volume of $FR(u)$ is at most $rn$. We can then fill $FR(u)$, and therefore $z$, in $H$ with volume at most $\delta^{(k-1)}_H(rn)$. 

\end{proof}

Note that under the area distortion definition used by Gersten, the upper bound is instead $AD \preceq n\delta_H$: this happens because the frontier may be disconnected, which introduces the presence of a summation. However, it is conjectured (see, for example, \cite{annals1}) that Dehn functions are super-additive, in which case we would regain the simpler bound $AD \preceq \delta_H$. When we choose homology, however, this complication disappears, because there is no requirement that the boundary be connected.\\

\begin{theorem}
\label{dehnlowerbound}
Let $H \subset G$ be $F_k$ groups, and let $\delta_H$ and $\delta_G$ be their respective $(k-1)$-order Dehn functions. Suppose $\delta_G$ is an invertible function. Then $\dist^k_{(G, H)} \geq \delta_H \circ \delta^{-1}_G.$
\end{theorem}

\begin{proof}
Let $z_n$ be a sequence of $(k-1)$-cycles with $$V_H(z) = \delta^{-1}_G(n)$$ and $$FV_H(z) = \delta_H(\delta^{-1}_G(n)),$$ that is, a sequence of maximally ``hard to fill'' cycles.  Then we know that $$FV_G(z) \leq \delta^{-1}_G(\delta_G(n)) = n.$$
\end{proof}

\subsection{Subgroups}

\begin{theorem}
\label{subsub}
Suppose $K \subset H \subset G$ are $F_k$ groups. Then:

\begin{itemize}
\item[(i)] $\dist^k_{(G, K)} \preceq \dist^k_{(H, K)} \circ \dist^k_{(G, H)}$

\item[(ii)] $\dist_{(H, K)}^k \preceq \dist_{(G, K)}^k$

\end{itemize}
\end{theorem}

\begin{proof}
For $(i)$, let $z$ be a $(k-1)$-cycle in $K$ so that $FV_G(z) \leq n$. Then $FV^k_H(z) \leq \dist^k_{(G, H)}(n)$ and so $FV^k_K(z) \leq \dist^k_{(K, H)}(\dist^k_{(G, H)}(n))$.

For $(ii)$, once again let $z$ be a $(k-1)$-cycle in $K$, but now suppose that $FV_H(z) \leq n$ but $FV_K(z) = \dist^k_{(K, H)}(n)$. Since $H \subset G$, a filling in $H$ is also a filling in $G$, so $FV_G(z) \leq n$ as well. Thus we've constructed examples of cycles in $K$ whose filling volume in $G$ is at most $n$, but whose filling volume in $K$ is $\dist^k_{(K, H)}(n)$, so this is a lower bound for the distortion of $K$ in $G$.

\end{proof}

This theorem is of particular interest if one of the embeddings is undistorted. 

\begin{cor}\label{undistsub}
Given $K \subset H \subset G$ as above, 

\begin{itemize}
\item{(i)} if $\dist^k_{(G, H)}$ is linear, then $\dist^k_{(H, K)} \asymp \dist^k_{(G, K)} $.
\item{(ii)} if $\dist^k_{(H, K)}$ is linear, then $\dist^k_{(G, K)} \preceq \dist^k_{(G, H)}$.
\item{(iii)} if $\dist^k_{(G, K)}$ is linear, then $\dist^k_{(H, K)}$ is linear.
\end{itemize}
\end{cor}

For an example applying these theorems, see \S \ref{heis}.

\subsection{Computability of distortion functions}

Papasoglu shows in Proposition 2.3 of \cite{papasoglu00} that area distortion functions are always computable. This contrasts significantly with the length distortion and first-order Dehn function cases, in which uncomputable functions can be obtained (see \cite{farb94}). The reason for this is that we start with objects we already know are the boundary of some chain; thus through brute force we will eventually be able to produce a filling.

In higher dimensions, Papasoglu's theorem generalizes:

\begin{theorem}
Given $F_p$ groups $H \subset G$, the function $\dist_{(G,H)}^k$ is computable for $k \geq 2$.
\end{theorem}

\begin{proof}
We can make a list of all ways of combining at most $n$ $k$-cells to make chains in $X_G$; now pick the subset whose boundary lies entirely in $X_H$. Because we know these must be trivial in $X_H$, we can find a volume (and thus a minimal volume) of a filling in $X_H$. 
\end{proof}

While this means that volume distortion is in some ways more nicely behaved than length distortion, \cite{madlenerotto85} provides the following theorem:

\begin{theorem}[Madlener-Otto]
Given a computable function $f$, there exists an example of a subgroup, group pair with area distortion bounded below by $f$.
\end{theorem}

\section{Examples}
\label{exs}

\subsection{Motivational examples}

There are many geometric examples which demonstrate the concept of area distortion; we will cover some of them here.

\subsubsection{Hyperbolic space}
Consider a horosphere in three-dimensional hyperbolic space, $\hyp^3$. Hyperbolic space has a linear isoperimetric function, so that a loop of length $l$ can be filled in area approximately $l$. When restricted to the horosphere, however, we encounter Euclidean geometry, for which we have quadratic isoperimetric function. As a result, we can find loops with area $l$ in $\hyp^3$ but $l^2$ in the horosphere, giving quadratic area distortion. We cannot hope for larger distortion, since the quadratic isoperimetric function provides an upper bound as well.

\subsubsection{Sol geometry}\label{sol}

Cosider the 3-dimensional Riemannian manifold $Sol$  that is topologically $\R^3$ and has the metric $ds^2 = \lambda^{-2t}dx^2 + (\frac{1}{\lambda})^{-2t}dy^2 + dt^2$. If we project $y$ to zero, we obtain a hyperbolic plane; projecting $x$ to zero gives an upside-down hyperbolic plane, and projecting $t$ to zero gives a Euclidean plane. Then the $xy$-plane is exponentially length-distorted in Sol, since the point on the $x$-axis of distance $n$ from the origin can be reached via a geodesic in the hyperbolic plane of the $xz$-axis of length approximately $\log{n}$. The $y$-axis is similarly distorted, but the geodesics now travel down instead of up. However, when we consider area, the two factors cancel each other out: it is exactly as good in the $x$-direction to go up as it is bad in the $y$-direction. Given any chain $z$, we know we can find the volume of $z$ as
\begin{align} 
V(z) &= \int \sqrt{(\lambda^{-t}(\frac{1}{\lambda})^{-t}dxdy)^2 + (\lambda^{-t}dxdt)^2 + ((\frac{1}{\lambda})^{-t}dydt)^2 }\\
&= \int \sqrt{(dxdy)^2 + (\lambda^{-t}dxdt)^2 + ((\frac{1}{\lambda})^{-t}dydt)^2 }
\end{align} 

Note that if we project $t$ to zero, we lose the second and third term, which can only decrease the overall volume; the $dxdy$ term is unaffected since the scalars cancelled out. Thus projecting $t$ to zero can only decrease area, which means that the $xy$-plane is area undistorted in $Sol$.

\begin{figure}
\begin{center}
\includegraphics[scale=.4]{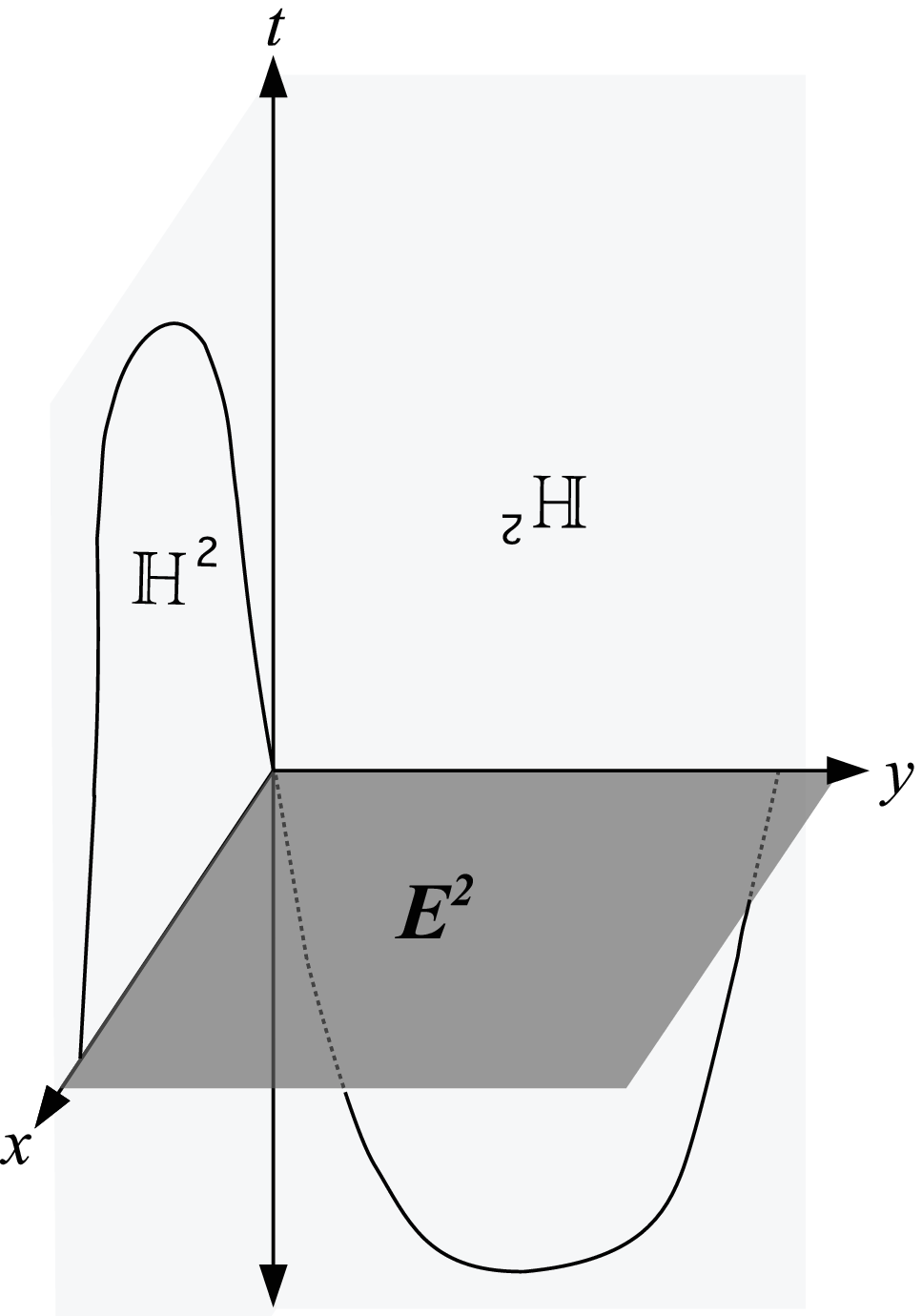}
\hspace{2cm}
\includegraphics[scale=.4]{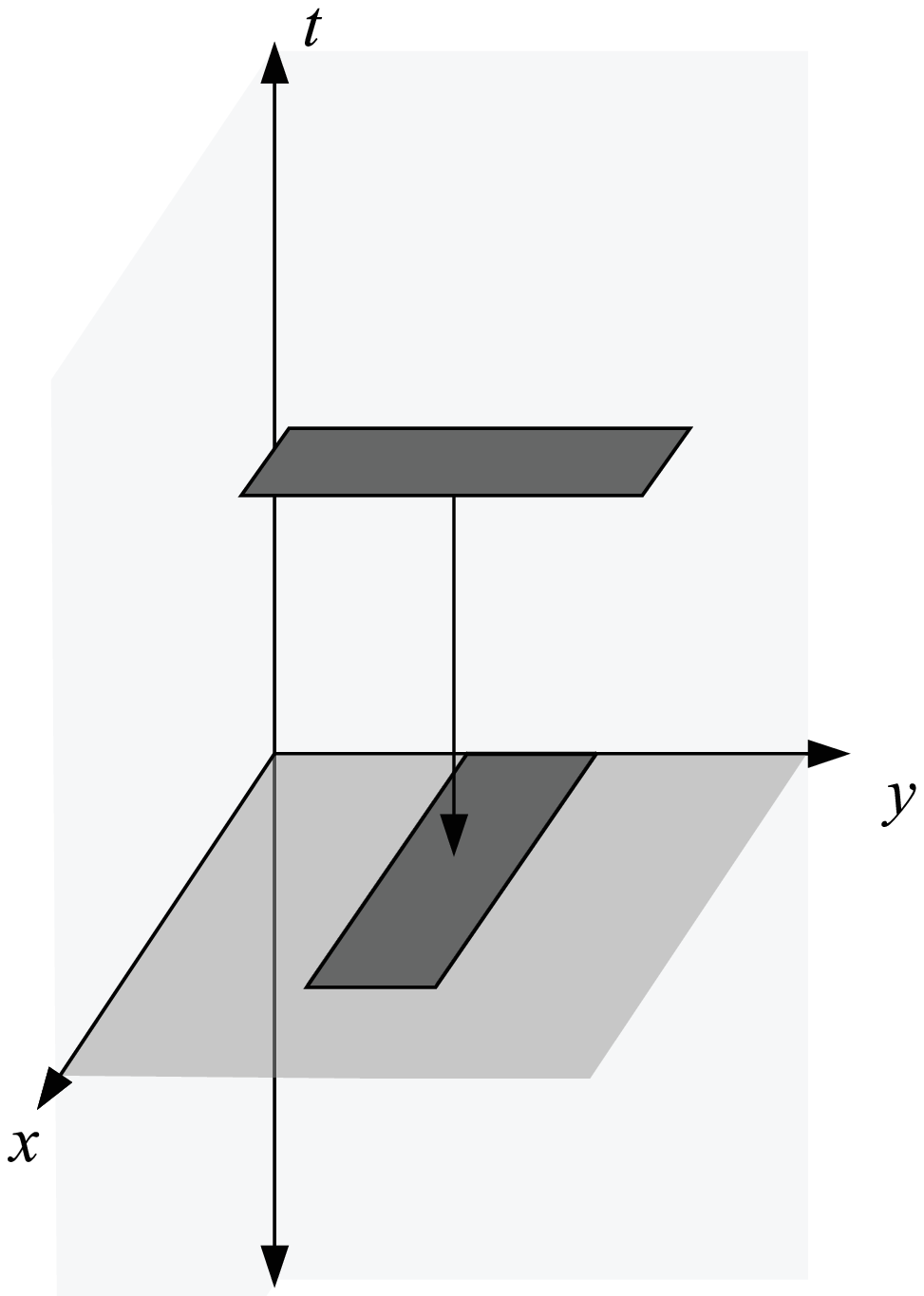}
\end{center}
\caption{In $Sol$, length is distorted (left) but area is not.}
\end{figure}

\subsubsection{Sullivan's theorem}\label{sull}
This has a generalization due to Sullivan. If $M$ is a 3-manifold and $\mathcal{F}$ is a codimension 1
foliation on $M$ which is transversely oriented, and such that there
is a transverse closed curve through every leaf, then there exists a Riemannian metric on $M$ for which every leaf of
$\mathcal{F}$ is quasi-area minimizing. As a special case, $\Z^2$ in $\Z^2 \rtimes_\phi \Z$ has undistorted area when $\phi
\in GL(n, \Z)$. Gersten proves this using a concept
he calls complexity; we will adjust the definition slightly, and
generalize it to higher dimensions.

\subsection{Complexity}
\label{cxitysection}

The concept of \textit{complexity} is defined in \cite{gersten96} for area distortion of $G$ in $G \rtimes_\phi \Z$, where $\phi$ is an automorphism of $G$ and $G$ is finitely presented. 

Let $G$ be $F_k$ and let $\phi$ an be automorphism of $G$. Then $\phi$ can be thought of as a map on the edges of a CW-complex $Y = K(G,1)$, where an edge labeled $s$ is sent to a word representing $\phi(s)$. Then each 2-cell is sent to a closed loop, and so we can fill it in some way in $Y$. Choose one of minimal area for each 2-cell, and call this $\phi(r)$. Continue with this process, inductively extending the map to the $n$-skeleton of $Y$ given the image of the $(n-1)$-skeleton, until we have a map $\phi: Y \to Y$. We can then lift this to a map $\tilde{\phi}: \widetilde{Y} \to \widetilde{Y}$.

\begin{figure}
\begin{center}
\includegraphics[scale=.4]{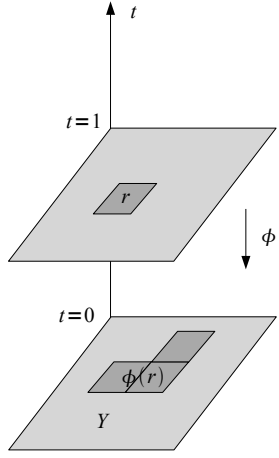}
\end{center}
\caption[Projection via $\phi$]{We can project to height zero via iterations of $\phi$.}\label{height}
\end{figure}

We can use this map to construct a $K(G \rtimes_\phi \Z, 1)$, given by $$X = Y \times I/[(y,0) \equiv (\phi(y), 1)],$$ which has universal cover that is setwise given by $\tilde{X} = \tilde{Y} \times \R$, made by taking a copy of $\tilde{Y} \times I$ for each integer, and identifying $(y,1)_i$ with $(\phi(y), 0)_{i+1}$. This construction gives us a natural projection $\pi: \tilde{X} \to \R$, where we call $\pi(x)$ the \textit{height} of the point $x \in \tilde{X}$. We can describe cells in this complex as follows: either they are cells inherited from $\tilde{Y}$, at a height $h$, or they are built inductively, with two-cells having boundary $tst^{-1}\phi(s)$ and higher-dimensional cells made by constructing the cells on the boundary and then filling them. See Figure~\ref{height} for a pictoral representation.

Let $c(\phi)$ be the maximal $k$-volume of the image of a $k$-cell under $\tilde{\phi}$. Since $G$ is $F_k$, $c(\phi)$ is finite. We will call $c(\phi)$ the \textit{$k$-complexity} of $\phi$.

\begin{remark}
Gersten defines (2-)complexity in a somewhat different manner: in \cite{gersten96}, the complexity of a map is the sum of the volumes of the images of all $2$-cells, minus the number of $2$-cells. His definition of complexity zero will coincide with our definition of complexity one. With his definition, one can only relate the distortion to complexity when the complexity is zero. We will be able to create a more general upper bound on the distortion, which depends on the value of $c(\phi)$.

\end{remark}

\noindent {\bf Theorem \ref{complexity}.}
Let $G$ be an $F_k$ group and let $\phi$ be an automorphism on $G$, and define $m = \max \{c_k(\phi), c_k(\phi^{-1})\}$, with $c_k(\phi)$ as defined above. Then the distortion of $G$ in $G \rtimes_\phi \Z$ is bounded above by $n \cdot m^n$.

\vspace{.5cm}

\noindent {\bf Corollary \ref{cxitycor}.}
When $\phi$ has complexity $m = 1$, then $G$ is $k$-volume undistorted in $G \rtimes_\phi \Z$.

In particular, if a $K(G,1)$ has only one $k$-cell, then $G$ is $k$-volume undistorted in $G \rtimes_\phi \Z$ for any automorphism $\phi$.

\begin{proof}[Proof of Corollary \ref{cxitycor}]
The first statement is trivial; the second is proven by noting that $\phi$ must send this $k$-cell to itself, because it induces an automorphism on on the $k$-homotopy of the $(k-1)$-skeleton of the space. Thus the $k$-complexity is 1.
\end{proof}

\begin{remark}

Theorem B of \cite{gersten96} proves the corollary in the case of $k = 2$ an Theorem 5.1 of \cite{gersten96} is related to Theorem \ref{complexity}; however, Gersten formulates his bound in such a way that Corollary \ref{cxitycor} does not follow from Theorem \ref{complexity}. 
\end{remark}

\begin{proof}[Proof of Theorem \ref{complexity}]

Let $G$ and $\phi$ be as in the theorem, and let $\tilde{X}$ be the universal cover of a $K(G \rtimes_\phi \Z, 1)$ as constructed above, with the height projection $\pi: \tilde{X} \to \R$. Let $z$ be a $(k-1)$-cycle in $\tilde{Y}$ with $u$ a $k$-chain in $\tilde{X}$, such that $\partial u = z$, and let $n = V^k(u)$. Assume $u$ has no connected components which are cycles, since these could be removed to decrease the volume of $u$, still giving a cycle with boundary $z$. Denote by $\triangle_p$ the subset of $u$ at height $p$, that is, $\triangle_p = \pi^{-1}(p) \cap u$. 

Since $u$ is closed, the image of $u$ under $\pi$ is a bounded subset of $\R$. Note that $\tilde{Y} = \pi^{-1}(0)$, so if $\pi(u) = 0$, we are done because $u$ is actually a chain in $\tilde{Y}$. Otherwise, we wish to use $u$ to construct a new chain that is in $\tilde{Y}$, whose volume is bounded above by $nm^{n}$. To do this, we will first consider $\pi^{-1}((0, \infty))$, and then $\pi^{-1}((-\infty, 0))$, which will work similarly.

If $\pi^{-1}((0, \infty))$ is empty, we proceed directly to the second set. Otherwise, we may choose some $p \in \R^{+}$ such that $p$ is not an integer, $\triangle_p \neq \emptyset$, and $\triangle_{p+1} = \emptyset$, that is, a height near the top of $u$. Let $U$ be the set of all cells in $u$ which intersect $\triangle_p$. This set may consist of a number of different connected components, but we need to break it up a bit more carefully: partition $U$ into subsets $U_1, U_2, \cdots, U_l$ so that in each subset, any two cells can be connected to each other in $u$ without going below height $p$. Then any two cells in the same connected component of $U$ will be in the same partition, but a partition may contain more than one connected component in $U$. 

Each $U_i$ then separates $u$ into two pieces, and gives a homology between some $(k-1)$-cycle $v_i$ at height $h = \lceil p \rceil$ and the $(k-1)$-cycle $\phi(v_i)$ at height $h-1$. The cycle $v_i$ must be the boundary for some sub-chain $c_i$ of $u$; by our choice of $p$, $c_i$ must be entirely at height $h$ and thus lie in a copy of $\tilde{Y}$. Then we can perform surgery on $u$, removing $U_i$ and $c_i$ and replacing them with $\tilde{\phi}(c_i)$ at height $h-1$. After doing this for each $i$, we have created a new $k$-chain $u'$ with height one less than the height of $u$. We may continue to do this until the maximal height of our new chain is zero.

We then do the same thing with $p < 0$ with the obvious adjustment on the choice of height, where we now use $\tilde{\phi}^{-1}$ to move the chain upward.

We must now calculate how much this surgery has increased the area. Any $k$-cell that had been at height $h$ has now been moved to height $0$, each time multiplying the volume by at most $m$, giving a total volume of $m^h$. How big can $h$ be? At most $n$, since there must be at least one cell at each height for a top cell to be connected to $z$. Thus our $n$ k-cells have been replaced by at most $nm^n$ cells at height zero. 

\end{proof}

Note that this upper bound is often much larger than the actual distortion. By Theorem~\ref{dehnupperbound}, if $m > 1$, then the $k^{th}$-order Dehn function of the subgroup must be greater than exponential for the complexity bound to be greater than the one provided by the Dehn function. One problem is the height---the upper bound of $n$ is almost certainly too large. In \S \ref{abc}, we will find other ways to bound the height in a particular class of examples so that $m^h$ can be made much smaller.

\subsection{Heisenberg groups}
\label{heis}
One special case of a group $G = \Z^2 \rtimes_{\phi} \Z$ is the
{ \em Heisenberg group}, where $\phi = \left(\begin{array}{cc}1 & 0 \\ 1 & 1\end{array}\right).$ This is also commonly written as 
$$\Heis^3 = \la x, y, z \mid [x, y] = z, [x, z] = [y, z] = 1 \ra.$$ 
(Note that the dimension ``acting'' on the $\Z^2$, usually denoted $t$ above, is $x$ here.) We can create ``higher-dimensional''
Heisenberg groups $\Heis^{2n + 1}$ with pairs of generators $x_i$,
$y_i$, along with $z$, such that each commutator $[x_i, y_i]$ is $z$,
and all other pairs commute. Note that any Heisenberg group is embedded in Heisenberg groups of higher order.

These groups are interesting to us in part because, while $\Heis^3$ has cubic Dehn function (see \cite{echlpt}), all higher dimensional Heisenberg groups have quadratic Dehn function, a theorem proven analytically by Allcock in \cite{allcock98} and later combinatorially by Ol'Shanskii and Sapir in \cite{olshanskiisapir}. Thus, by Theorem~\ref{dehnlowerbound}, the distortion of $\Heis^3$ in $\Heis^5$ or any higher-order Heisenberg group is at least $n^{3/2}$.

Further, we can construct an upper bound as follows: consider the intermediate group $\Heis^3 \times \Z$, where $\Z$ is generated by $y_2$. Then we have
$$\Heis^3 \subset \Heis^3 \times \Z \subset \Heis^5 = (\Heis^3 \times \Z) \rtimes \Z,$$
 with the action given by $x_2$ commuting with $\Heis^3$ and sending $y_2$ to $y_2z$. The first containment is undistorted, so by Theorem~\ref{undistsub} (ii), the distortion of $\Heis^3$ in $\Heis^5$ is at most the distortion of  $\Heis^3 \times \Z$ in $\Heis^5$. Directly applying Theorem~\ref{complexity} would give us an exponential upper bound, but we can modify it a bit: notice that the
automorphism preserves all relators except commutators $[y_2, s]$ for $s = x_1, y_1$, each of which goes to two relators: a copy of itself, and the commutator $[z,y_2]$. Thus repeated applications of the automorphism only increase the image by one. Thus the final volume is at most $nh$, for height $h$, which is at most
$n$. Therefore $n^2$ is an upper bound for the area distortion, an improvement over the bound of $n^3$ given by the Dehn function of $\Heis^3$.

\begin{conjecture}\label{heisconj}
The area distortion of $\Heis^{3}$ in $\Heis^5$ is $n^{3/2}$.
\end{conjecture}

The reason for this is that the upper bound fails to take into account the ``side area'' coming from any filling; conceptually, any $y_2$ edge of height $h$ ought to be creating a side with $(h-t)$ relators at height $t$, creating a total area of $h^2$.

\subsection{Abelian-by-cyclic groups}
\label{abc}

Note that in this section, upper bounds are found without regard to the topology of the objects, and lower bounds are given by filling spheres with balls, so that the methods and results described work equally well if distortion is defined via homotopy rather than homology.

A group $\Gamma$ is \textit{abelian-by-cyclic} if there is an exact sequence $1 \to A \to \Gamma \to \Z \to 1$. By a theorem of Bieri and Strebel, given a finitely presented, torsion-free abelian-by-cyclic group, there is an $m \times m$ matrix $M$ with integer entries so that $\Gamma$ has the presentation 
$$\Gamma_M = \la x_1, \cdots x_m, t  \mid [x_i, x_j] = 1, tx_it^{-1} = \phi(x_i) \text{ for } 1 \leq i, j \leq m \ra,$$ 
where $\phi$ is a homomorphism taking $x_i$ to $x_1^{a_1}x_2^{a_2} \cdots x_m^{a_m}$, where the $a_j$ form the $i$th column of $M$. Several ideas from  \cite{farbmosher00} will help us to find the $k$-volume distortion of $\Z^m$ in $\Gamma_M$. 

We can construct a space $X_M$ on which $\Gamma_M$ acts properly discontinuously and cocompactly by isometries, so that $X_M$ and $\Gamma_M$ are quasi-isometric. Topologically, this space is $\R^m \times T_M$, where $T_M$ is a directed tree with one edge entering each vertex and $det(M)$ edges leaving each vertex. To give $X_M$ a metric, assign $T_M$ a metric with each edge having length 1,  and fix a particular vertex $v_0$ of $T_M$. This choice of $v_0$ gives us a height function from $T_M$ to $\R$, which we can extend to a height function $h: X_M \to \R$, where $h(v_0) = 0$. Then $(\Z^m, \Z^m \rtimes \Z)$ and $(R^m \times \{v_0\}, X_M)$ are quasi-isometric as pairs.

We can also consider the continuous Lie group $G_M = \R^m \rtimes_M \R$. Here multiplication is given by $(x, t) \cdot (y, s) = (x + M^ty, t+s)$. This space is a Riemannian manifold with left-invariant metric 
$$g_{ij}(x,t) = \left(\begin{array}{cc}(M^{-t})^TM^{-t} & 0 \\0 & 1\end{array}\right).$$

While this metric involves $M$, choosing any power of $M$ will give us a quasi-isometric space, so that we may replace $M$ with $M^2$ (for example, if $det(M) < 0$) or $M^{-1}$ (if $|det(M)| < 1$; this amounts to flipping the space vertically) as we wish.

We will then be able to integrate using this metric to find volumes. We also have a natural height function on $G_M$ given by the last coordinate. 

\begin{figure}
$$\xymatrix{ & X_M \ar[ld]_{g_M} \ar[dd]^h \ar[rd]^{\pi_M} & \\
G_M \ar[rd]_{\pi} & & T_M \ar[ld] \\
& \R & }$$
\caption{Relationship between $X_M$ and $\Gamma_M$}\label{gmxm}
\end{figure}
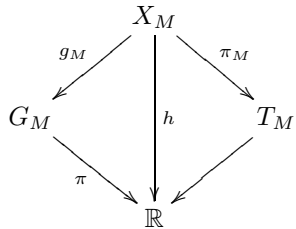

While $X_M$ and $G_M$ generally differ---if $\det(M) > 1$, then $X_M$ is not a manifold---there is a relationship between them that allows us to say that they have the same $k$-volume distortion. This relationship can be seen in the form of a commutative diagram, as seen in Figure \ref{gmxm}.

\begin{lemma}\label{GXeq}
The volume distortion function $\dist_X$ of $\R^m \times \{v_0\}$ in $X_M$ is equivalent to $\dist_G$ of $G_0 = \pi^{-1}(0)$ in $G_M$.
\end{lemma}

\begin{proof}
Choose some cross-section $\hat{g}: G_M \to X_M$ so that the image contains $v_0$; this gives an isomorphic embedding of $G_M$ into $X_M$. Given a chain in $X_M$ with boundary in $\R^m \times \{v_0\}$, we can map it under $g_M$ to a chain in $G_M$ with boundary in $G_0$. We can then fill the chain in $G_0$ and use $\hat{g}$ to pull this back to a filling in $\R^m \times \{v_0\}$; this tells us $\dist_X \preceq \dist_G$. The same process can be used to show the inequality holds the other way, giving us equivalent functions. 
\end{proof}

Note that the pulling back via $\hat{g}$ is possible only because our filling lies entirely at height zero, on which $g$ and $\hat{g}$ act as inverses. (This issue is important when considering Dehn functions in these groups, where the boundary is no longer restricted to a particular vertex in the tree; see \cite{bradyforester}.)

We may now simplify the situation to considering the height zero subspace $\R^m$  inside of $G_M$. By \cite{farbmosher00}, this group is quasi-isometric to $G_N$, where $N$ is the {\em absolute Jordan form} for $M$, that is, a matrix with the absolute values of eigenvalues along the diagonal and ones and zeroes elsewhere, in accordance with the Jordan form. This quasi-isometry preserves the height-zero subspace, so we may restrict our attention to matrices in Jordan form with positive real eigenvalues.

\subsubsection{Diagonalizable matrices}

Suppose $M$ is in absolute Jordan form, $det(M) \geq 1$. Call the $(i,i)^{th}$ entry $\lambda_i$. If $M$ has ones on the superdiagonal, the situation gets somewhat more complicated; we shall first restrict our attention to the case that $M$ is in fact diagonal. In this case, the geometry of the resulting Riemannian metric is particularly easy to understand. Topologically, we have a space of the form $\R^m \times \R$, where the last coordinate, denoted $t$, will be considered the height. The metric either expands (if $\lambda_i < 1$) or contracts (if $\lambda_i > 1$) the $x_i$ direction as the height increases. We have the metric:
$$ds^2 = dt^2 + \sum_{i = 1}^{m} \lambda_i^{-2t} dx_i^2.$$

Further, given a map $g: \sigma^k \to G_M$, with image given by $(g_1, g_2, ..., t)$, the $k$-volume of $g$ in $G_M$ is
$$\int_{\sigma^k} (\sum \lambda_I^{-2t} |D_Ig_x|^2)^{1/2} dx,$$
where $I$ is a choice of $k$ of the basis vectors, $|D_Ig_x|$ is the determinant of $Dg$ restricted to those $k$ vectors, $\lambda_I$ is the product of the $\lambda_{i}$ of $M$ for $i \in I$. Note that the choice $m+1$ gives the $t$ direction; thus $\lambda_{m+1} = 1$.

\begin{theorem}
\label{abcdist}
Let $M$ be a diagonal $k \times k$ matrix, with $(i,i)^{th}$ entry $\lambda_i \in \R^+$ and determinant $d > 1$ and at least two eigenvalues off the unit circle;.  Then the $k$-volume distortion of the height-zero copy of $\R^k$ in $G_M$ is the function
$$ \dist^{(k)}(n) = n^{1 + \frac{\log{d}}{log{\alpha}}}, \textrm{ where }  \alpha = {(\prod_{i = 1}^k \max \{d, \lambda_i\})}/d.$$
\end{theorem}

This proof has benefitted, both in scope and simplicity, from ideas provided by Brady and Forester in \cite{bradyforester}.

\begin{proof}

To simplify later calculations, define $p_i = d/\lambda_i$ and $p = \prod_{i = 1}^{k} \min\{p_i, 1\}$. We can then compute that $\alpha$ can also be written as 
$$\alpha = \frac{d^{k-1}}{p};$$
this formulation requires more notation, but will better match the approach taken in the proof.

We begin by showing this function is an upper bound. Suppose $u$ is a $k$-chain in $G_M$ with volume $n$, and boundary contained in $\R^k$. Then it will suffice to prove that we can fill $\partial u$ in $\R^k$ with volume $n^{1 + \log{d}/\log{\alpha}}$.

First, break $z$ into two pieces: the ``low'' piece $u_L$ consisting of all of $u$ with height less than $h = \log(n)/\log(\alpha)$, and the ``high'' piece $u_H$. 

The projection $\pi_t$ sending $t$ to 0 sends $u$ to a filling of $\partial u$, but increases the volume. We know that the volume of $\pi_t(u_L)$ is an increase of the volume of $u_L$ by a factor of at most $d^{h} = n^{\log{d}/\log{\alpha}}$, giving us a total volume as desired. Now we need only bound the total volume of $w$, the regions interior to $\partial u$ covered only by $u_H$. The process for doing this is illustrated in Figures \ref{projpf0}--\ref{projpf5}.

\begin{figure}
\begin{center}
\includegraphics[scale=.3]{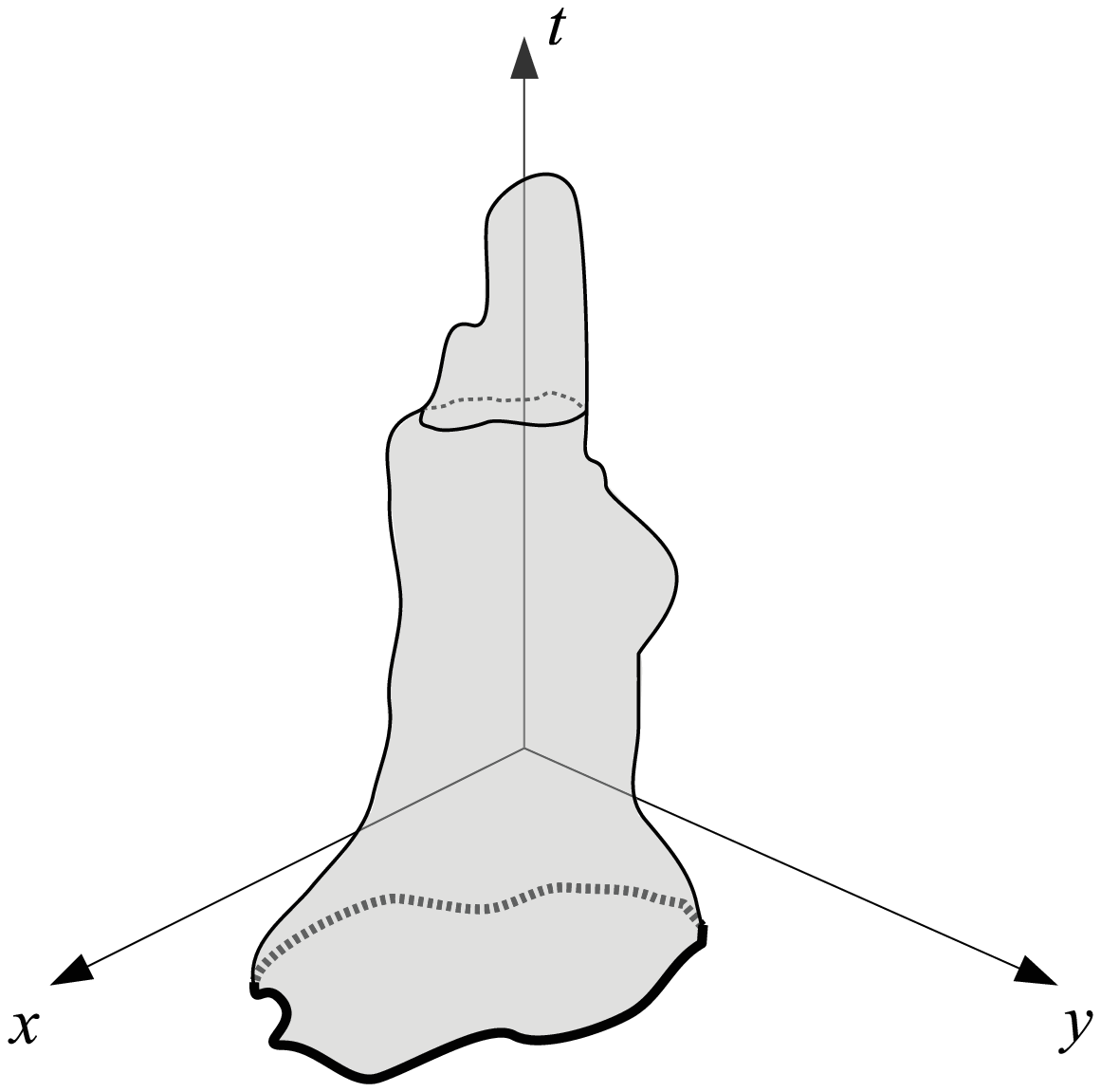}
\hspace{1cm}
\includegraphics[scale=.3]{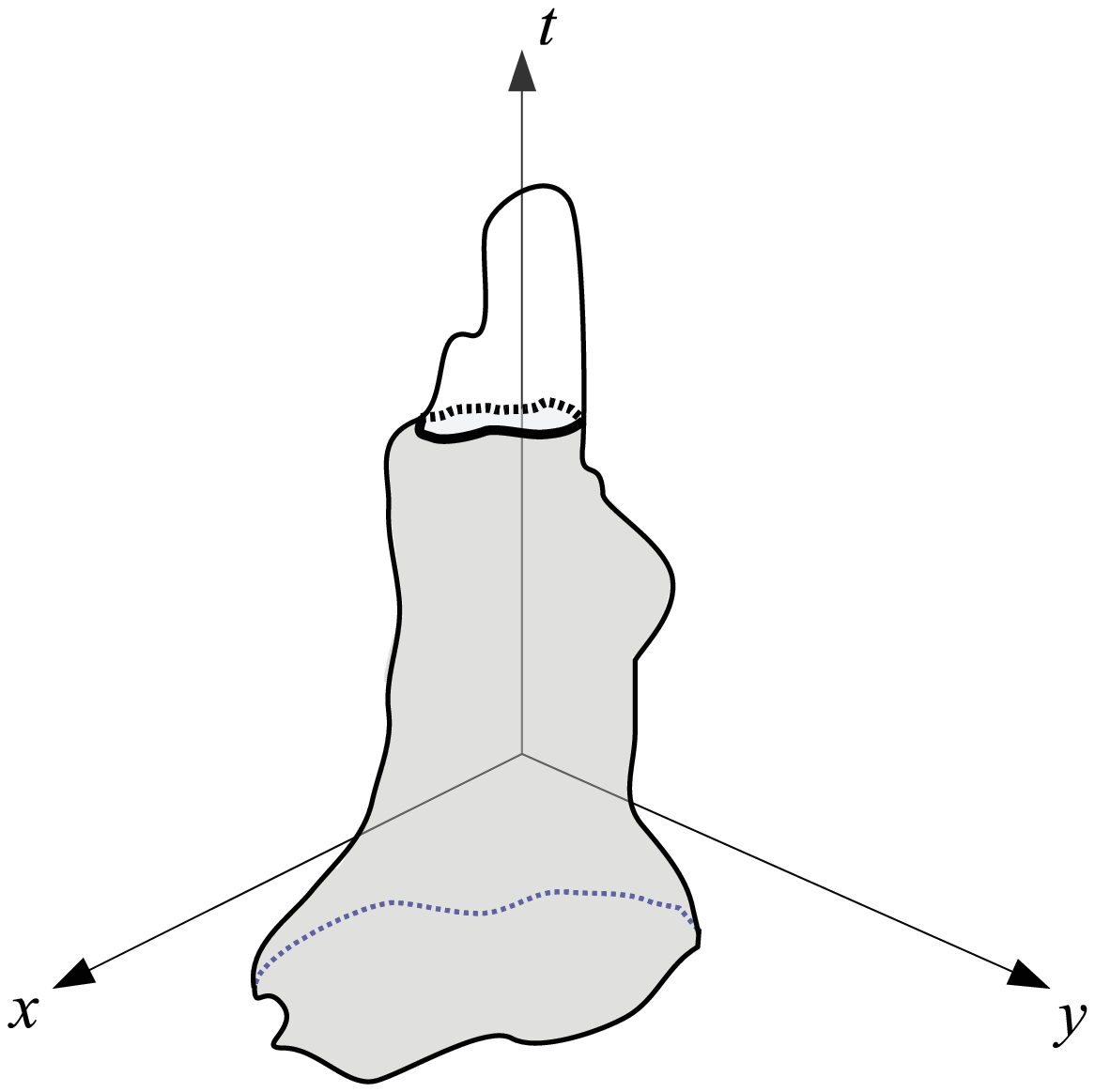}
\end{center}
\caption{Start with a cycle with a filling of volume $n$ and cut off at height $h$.}\label{projpf0}
\end{figure}

\begin{figure}
\begin{center}
\includegraphics[scale=.3]{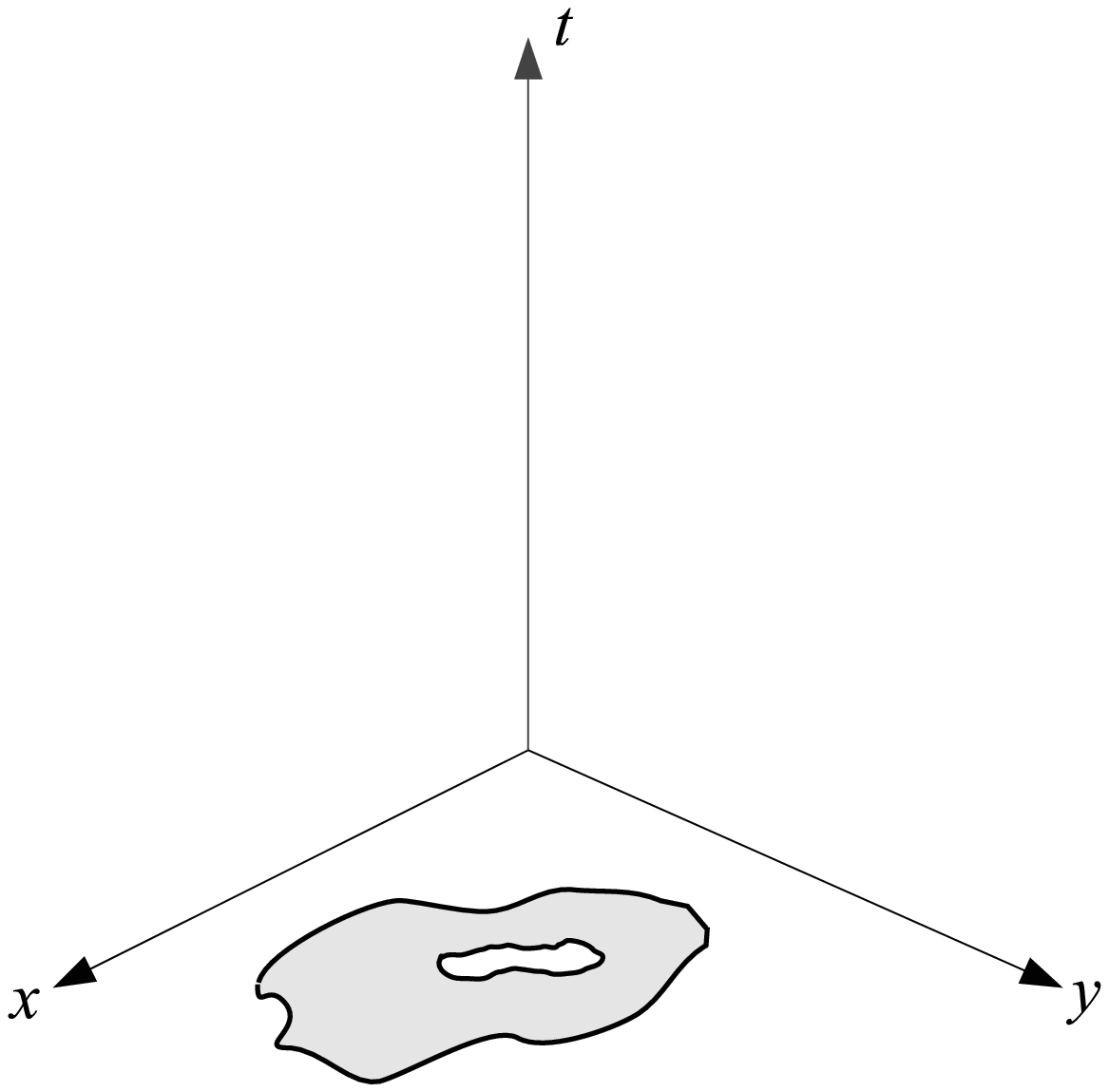}
\hspace{1cm}
\includegraphics[scale=.3]{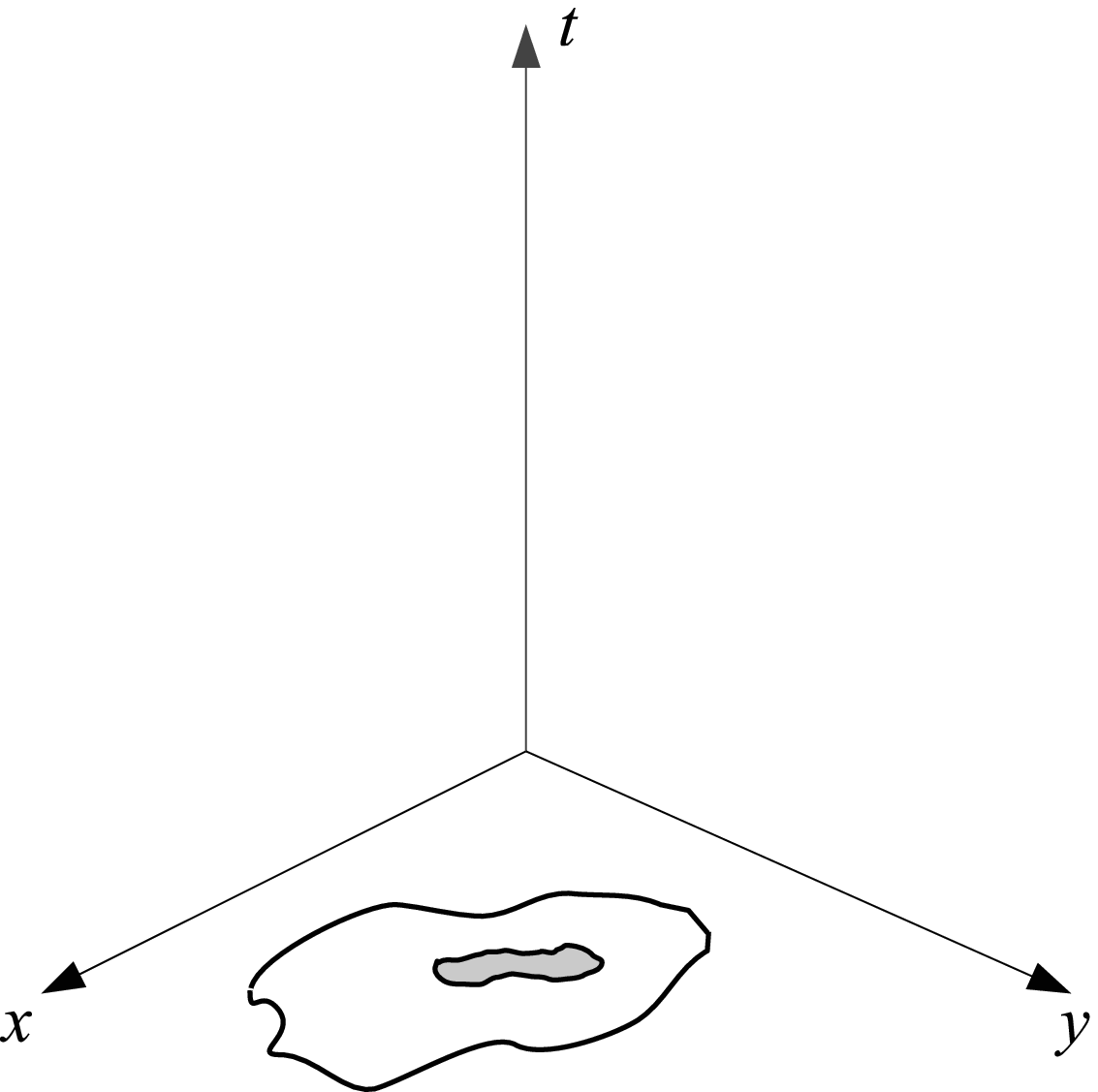}
\end{center}
\caption{Project low piece down; we want to bound the volume of the pieces not filled.}\label{projpf2}
\end{figure}

\begin{figure}
\begin{center}
\includegraphics[scale=.3]{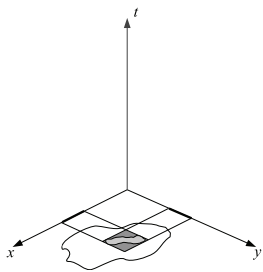}
\end{center}
\caption{Do this by projecting each coordinate to zero; bound the volumes of the projections.}\label{projpf4}
\end{figure}

\begin{figure}
\begin{center}
\includegraphics[scale=.3]{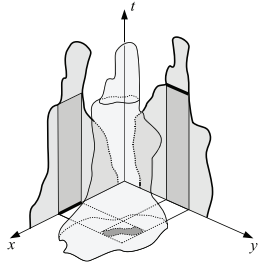}
\end{center}
\caption{The projection must contain a cylinder; this bounds the volume at height zero.}\label{projpf5}
\end{figure}

\begin{figure}
\begin{center}
\includegraphics[scale=.3]{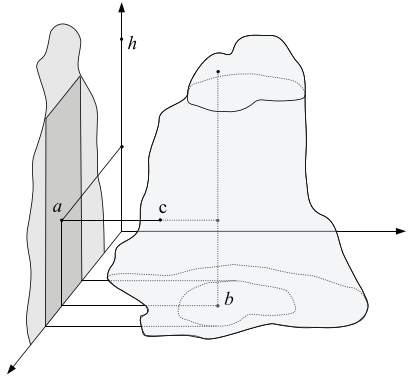}
\end{center}
\caption{A cylinder with base $w_i$ (and thus base volume $v_i$) and height $h$ must appear in the projection.}\label{cylinder}
\end{figure}

First note that the map $\pi_i$ projecting $x_i$ to 0 is volume non-increasing. We will finds bounds $V_i$ on the $(k-1)$-volume of $w_i = \pi_i(w)$; this allows us to bound the final volume 

\begin{equation}\label{Vw}
V(w) \leq (\prod V_i)^{1/(k-1)}
\end{equation}

In particular, we will show that if $p_i > 1$, then $V_i \leq C_in$, where $C_i$ is independent of $n$, and if $p_i < 1$, then $V_i \leq C_inp_i^{h}$.  

Given these bounds on the $V_i$, notice that when we multiply them together, we get a factor of $p_i^h$ for each $p_i < 1$. This is exactly the definition of $p^h$. Thus
\begin{equation}\label{volumeproduct}
\prod V_i \leq Cn^kp^h 
\end{equation}

We can use the relationships between $p, \alpha, h$, and $n$ and basic properties of logarithms to find that 
\begin{equation}\label{ph}
p^h =  n^{\frac{(k-1)\log d}{\log \alpha} - 1}.
\end{equation}
We can now substitute (\ref{volumeproduct}) and (\ref{ph}) into (\ref{Vw}) to find
\begin{equation}\label{finalvol}
V(w) \leq (C'n^{(k-1)(1 + \frac{\log d}{\log \alpha})})^{\frac{1}{k-1}} = C'n^{1 + \frac{\log d}{\log \alpha}}.
\end{equation}

 This allows us to fill all of $\partial u$ at height zero with volume at most $(1+C)n^{1 + \log{d}/\log{\alpha}}$, proving our upper bound.

It remains to show that these bounds on the $V_i$ are valid. We will do so by showing that the projection $\pi_i(u)$ contains the cylinder with base $w_i$ and height $h$; the bound on $V_i$ then comes from the fact that the volume of the cylinder cannot be more than $n$. The process described for finding this cylinder is illustrated in the area distortion case in Figure \ref{cylinder}.

Without loss of generality, assume $i = k$. Let $a$ be a point in the cylinder, say with $(x, x_k,t)$-coordinates given by $(q,0,h_0)$. Then $(q,0,0)$ is interior to $u_k$, so there is some $q_k$ such that $b = (q,q_k,0) \in u$. Now consider the line $(q,q_k, t)$. This will intersect $w$ for the first time at some height $H > h$. In particular, the point $(q,q_i,h_0)$ is interior to $W$. Then the line $(q,x_i,h)$ must intersect $W$ at some point $c$. Thus $\pi_k(c) = a$.

The volume of the cylinder is given by $\int_0^h p_i^{-t}V_i dt$.  We will restrict the height further: if $p_i > 1$, we will consider only the cylinder from height 0 to 1; if $p_i < 1$ then we will consider the cylinder from height $h -1$ to $h$. In the former case we have:
\begin{displaymath}
\begin{aligned}
\int_0^1 p_i^{-t}V_i dt &= -\frac{1}{\log{p_i}}(p_i^{-t}V_i)|_o^1\\
 &= \frac{V_i}{\log{\pi_i}}(1 - p_i^{-1})
\end{aligned}
\end{displaymath}

Since this volume is less than $n$, we must have $$V_i \leq n (\log{p_i})(1-p_i^{-t})^{-1} = C_in.$$

Similarly, in the latter case we have:

\begin{displaymath}
\begin{aligned}
\int_{h-1}^{h} p_i^{-t}V_i dt &= -\frac{1}{\log{p_i}}(p_i^{-t}V_i)|_{h - 1}^{h}\\
&= \frac{V_i}{|\log{pi}|}p_i^{-h}(1 - p_i)
\end{aligned}
\end{displaymath}

Thus $V_i \leq |\log{p_i}| n p_i^{h}(1 - p_i)^{-1} = Cnp^{h}$, as desired.

\vspace{1cm}

To show this bound is sharp, we need to construct an example exhibiting this amount of distortion. We do this by costructing a $(k-1)$-dimensional box with side lengths chosen so that each projection $\pi_i$ gives an object with $(k-1)$-volume equal to the upper bound $V_i$ found above.

In order to do this, set $l_i = (\prod V_i)^{1/(k-1)}/V_j$. Now build a $(k-1)$-hyper-rectangle with the length in the $x_i$ direction equal to $l_i$.

This box can be filled in the ambient space by flowing the box up to the height $h$ and then filling the resulting box at height $h$. This height was chosen so that the volume obtained by flowing each side of the box is some $Cn$, and the the volume at height $h$ is also $Cn$. Thus the overall volume of this filling is bounded above by $C'n$.

The subspace is Euclidean, so we know the best filling for the box, which is exactly $\prod V_i = Cn^{1 + \frac{\log{d}}{\log \alpha}}$, as calculated in Equations \ref{Vw} through \ref{finalvol}, which is exactly the value of the volume distortion function.
\end{proof}

In the case $M$ has one eigenvalue off the unit circle, we must amend the bound somewhat. In particular, it is necessary to use the {\em Lambert $W$ function}, that is, the inverse of the function $e^nn$.

\begin{cor}\label{oneeig}
If $M$ is as above, but has exactly one eigenvalue off the unit circle, say $\lambda > 1$, then the $k$-volume distortion function is  $(\frac{n^{k}}{W(n)})^{1/(k-1)}$, where $W$ is the Lambert $W$ function. 
\end{cor}

\begin{proof}
The reason we cannot use exactly the proof of Theorem~\ref{abcdist} is that when we project $x_1$ to zero, we obtain a space with a Euclidean metric. Thus the volume of the cylinder will be $h \cdot V_1$. 

Because of this, we must change $h$ to $W(n)$ and set $V_1 = n/h$. Notice that all other $V_i$ will still be $n$, since these projections behave as before. With these changes, the argument for the proof of Theorem \ref{abc} works exactly, giving a $k$-volume distortion function of $\lambda^hn$, which is equivalent to $(\frac{n^{k}}{W(n)})^{1/(k-1)}$. 

\end{proof}

By Lemma \ref{GXeq} above, the distortion of $\Z^k$ in $\Gamma_M$ is the same as that of $\R^k$ in $G_N$, where $N$ is the absolute Jordan form for $M$, so we immediately have:\\

\noindent
{\bf Theorem \ref{gammamdist}.}
Let $M$ be an integer-entry $k$-by-$k$ diagonalizable matrix with
 $det(M) = d \geq 1$,and let $\lambda_i$ denote the absolute value of the $i^{th}$ eigenvalue. Then the $k$-volume distortion of $\Z^k$ in $\Gamma_M$ depends only on the eigenvalues of $M$. If $M$ has at least two eigenvalues off the unit circle, the volume distortion is 

$$\dist^{(k)}(n) \asymp n^{1 + \log{d}/\log{\alpha}}, \textrm{ where } \alpha = (\prod_{i = 1}^{k} \max \{\lambda_i, d\})/d.$$
  
If $M$ has exactly one eigenvalue off the unit circle, 
$$ \dist^{(k)}{n} \asymp (\frac{n^{k}}{W(n)})^{1/(k-1)}.$$

Otherwise, $\dist^{(k)}(n) \asymp n$.

\vspace{1cm}

The case in which $M$ has determinant $d = 1$ is covered by Theorem \ref{cxitycor}, since in this case $M$ gives an automorphism, and $\Z^k$ has a unique $k$-cell. This presents one extreme, the case in which volume is undistorted.

At the other extreme, when all eigenvalues are at least 1, and when at least two eigenvectors are greater than one, then the $k$-volume distortion is maximal, i.e. $n^{k/(k-1)}$, which is the $(k-1)$-order Dehn function for $\Z^k$.

We may wish to consider the $k$-volume distortion of $\Z^m$ in $\Gamma_M$ with $k < m$, that is, the distortion of a smaller-dimensional volume. This is bounded below by the largest $k$-volume distortion of $\Z^k$ inside the group we get by projecting the other $m-k$ dimensions to zero, as these projections are volume non-increasing.

\begin{cor}
For any intergers $m > 1$ and $1 < k < m$, there exists a pair $(G, H)$ with distorted $k$-volume but undistorted $m$-volume.
\end{cor}

\begin{proof}
Simply choose a group $G = \Gamma_M$ with $M$ an $m$-by-$m$ matrix with $det(M) = 1$ and at least one eigenvector off the unit circle, and let $H = \Z^m$.

\end{proof}

Ideally, we would like groups that exhibit stronger behavior: for example, a pair in which only the 3-volume is distorted. Examples exist for the area case, but have proven more difficult to construct in general.

\subsubsection{Other matrices}\label{nondiag}

When the matrix $M$ is not diagonalizable, the situation gets more complicated. In this case, the automorphism no longer preserves the eigendirections, but also changes lengths along other directions at a rate proportional to a polynomial in the height. This means that techniques involving projection become more difficult to use.

\begin{conjecture}
Given a matrix $M$ with at least one eigenvalue off the unit circle, the distortion of $\Z^k$ in $\Z^k \rtimes_M \Z$ is the same as that of $\R^k$ in $\R^k \rtimes_N \R$, where $N$ is a diagonal matrix with diagonal entries given by the norms of the eigenvalues of $M$.
\end{conjecture}

The idea behind this conjecture is that the exponential change in length created by the eigenvalues dominates the polynomial change given by the ones on the superdiagonal, and so, in a large-scale sense, we should be able to ignore the polynomial contribution.

This leaves one more case: when all of the eigenvalues are on the unit circle. In this case, there is no exponential growth coming from the eigenvalues, and so it is the polynomial effect that comes into play. 

Such a matrix will have absolute Jordan form in which each block has ones on the diagonal and superdiagonal, and zeroes everywhere else. Thus it can be described completely by the number of blocks, say $c$, and the size of each block, which we will denote by $a_i$ for $i \in \{1, 2, \dots, c\}$.

\begin{theorem}\label{blockdist}
Let $M$ be a matrix whose absolute Jordan form consists of $c$ blocks along the diagonal, in which each block is $a_i$-by-$a_i$ and consists of ones on the diagonal and superdiagonal and zeros elsehwere. Then the $k$-volume distortion of $\Z^m$ in $\Gamma_M$ is at least $n^{\beta/\alpha}$, where 

\begin{equation}
\begin{aligned}
\alpha &= (k - 1)\sum (a_i -1) + k, \text{ and} \\
\beta &= k\sum (a_i - 1) + k
\end{aligned}
\end{equation}

if $k \leq c$, and

\begin{displaymath}
\begin{aligned}
\alpha &= (k-1)(m-k) + 2k - c + \sum (b_i - k_i) \text{, and }\\
\beta &= k(m-k) + 2k - c + \sum (b_i - k_i)
\end{aligned}
\end{displaymath}

when $k > c$, where $k_i$ are chosen such that $k_i \leq b_i$, the $b_i$ are a subset of the $a_i$, and the choice of $k_i$ maximizes the value of $\alpha/\beta$.

Further, the $k$-volume distortion is bounded above by $n^{\frac{k}{k-1}}$.
\end{theorem}

\begin{proof}

The upper bound comes from the $k^{th}$-order Dehn function for Euclidean space.

To show the lower bound, we will construct examples with the desired distortion in the following general manner. We start with the boundary of a $k$-dimensional box, where each edge is in a basis direction $x_i$ with length $l_i$. We then apply the automorphism $M^h$ to this box to get an $(k-1)$-dimensional parallelepiped. We will embed this in $G_M$ at height 0; this will be the $(k-1)$-cycle that we consider. Call the map that gives us this cycle $f(u)$. It will often be beneficial to consider the edges of the box as vectors.

This cycle can be filled in at height $0$, and we can easily compute the volume of the resulting parallelepiped. On the other hand, we can think of $G_M$ as (topologically) $\R^m \times \R$ and fill the cycle with the function $(f(u), hz)$, together with filling the paralellepiped at height $h$. When $k=2$, this would look like filling the rectangle bounding the bottom of a box with its sides and top. 

The volume of the sides, which look like $\pi_i((f(u), hz))$, can be found by integrating the $(k-1)$-volume of the cross-section at a given height. We chose our chain so that this would be easy to compute: the volume of $f(u)$ at height $t$ in the space $G_M$ is the same as the Euclidean volume of $M^{-h}(f(u))$. This volume will be some polynomial in $h$ and the $l_i$; we will choose the values of the lengths $l_i$ and height $h$ so as to maximize the difference between the two volumes. To do this, we will set values so that the volumes of each of the sides and the top are all equal. Note that, while we do not prove that this maximizes the distortion, it still provides a lower bound for distortion, since we will have demonstrated an example with the distortion required.

It is convenient now to consider two cases: one, when $k \leq c$, and another when $k > c$. When two basis directions are chosen from the same Jordan block, the resulting vectors are not linearly independent, a fact that must be taken into account when computing volumes---it tends to significantly decrease the volume.

For this reason, if $k \leq c$, we will choose each vector from a different Jordan block. If a block has size $a$, then the image of the unit vector corresponding to the last column of this block is mapped under $M^h$ to a vector of the form $p_0(t)v_a + p_1(t)v_{a-1} + p_2(t)v_{a-2} + \cdots p_{a-1}(t)v_1$, where $v_i$ is the basis vector corresponding to the $i^{th}$ column of the block in question and $p_i$ is a polynomial of degree $i$. This will then be equivalent to a polynomial of degree $a-1$ in the height. 

If we choose our $k$ vectors to be the last vector in each of the $k$ largest Jordan blocks, the volume of the side where the $j^{th}$ side is sent to zero is given by 
$$ \int_0^h t^{\alpha_i} \frac{l_1l_2\cdots l_k}{l_j} = h^{\alpha_i + 1} \frac{l_1l_2\cdots l_k}{l_j}
$$
where
$$\alpha_i = \sum_{i \neq j} (a_i - 1).$$

These volumes will be equal when each $l_i = h^{\alpha_i}$. Under this choice, the volume of the filling in the ambient group is $h^{\alpha}$ where 
$$\alpha = \sum (\alpha_i + 1) = \frac{k-1}{\sum (a_i-1)} + k.$$

At height zero, the lengths of the vectors are now $h^{a_i - 1}l_i$, so the volume below is $h^{\sum a_i} \cdot l_1l_2 \cdots l_k = h^{\beta}$ where 
$$\beta = k\sum (a_i - 1) + k.$$ 

When we set the filling in the ambient group to be bounded above by $n$, then the area below, in terms of $n$, is $n^{\beta/\alpha}$.

Next we will find the volumes in the case $k > c$. We can, of course, also use this method to produce one bound for all choices of $k$ and $c$; however, if we know we are using all of the Jordan blocks, we can simplify the result somewhat to get a form that is less dependent on the sizes of the blocks than one might expect.

For the moment, let $a_i$ denote the sizes of the blocks from which only one vector is chosen, and $b_i$ denote the sizes of blocks in which multiple vectors are chosen, with $k_i$ representing the number chosen from this block. 

As before, we wish to find the $k$-volumes of the ``side'' pieces, which we do by integrating the Euclidean volume of $M^{-h}(f(u))$ from height $t = 0$ to height $t = h$. This volume is found by taking the determinant of the matrix whose columns are the vectors of the sides of the box (in each case, projecting one of the vectors to 0). When we take vectors from the same block in $M$, the linear dependence of the vectors significantly decreases the volume from what it would be if they were independent. Generally, if we are considering a side which contains the edges corresponding to the last $k_i$ vectors from a block of size $b_i$, the volume at height $t$ is $t^{b_i-k_i} \prod l_i$. Recall that one vector will be projected to zero; if it is one of these vectors, but not the last from the block, then the power of $t$ is increased by one, created by the ``gap'' in the dependence relations.

We now find the values of $l_j$ as powers of $h$ just as above: if the $j^{th}$ vector is from a block with only one chosen vector, then $l_j = h^{\alpha_j + 1}$ where
$$\alpha_j = \sum_{i \neq j} (a_i - 1) + \sum (b_i - k_i) +1$$

If $v_j$ is from a block with more than one vector, the power of $h$ is  
$$\sum (a_i - 1) + \sum(b_i -k_i) + q,$$

where $q$ is 1 when $v_j$ is the last vector in its block, and 2 otherwise.

With these values, the total volume of the filling is $h^\alpha$ where
$$\alpha = (k-1)\sum (a_i -1) + k\sum(b_i - k_i) + k + \sum (n_i - 1).$$

Notice that we can simplify some: we can think of 1 as the ``$k_i$'' for the blocks of size $a_i$. With this approach, the sum of the $k_i$'s is $k$. Further, the sum of the $a_i's$ and $b_i's$ is $m$, the dimension of the matrix. Then
\begin{displaymath}
\begin{aligned}
\alpha &= (k-1)[(\sum a_i + \sum b_i) - \sum k_i] + \sum (k_i - 1) + \sum(b_i - k_i)\\ 
&= (k-1)(m-k) +  2k - c + \sum(b_i - k_i)
\end{aligned}
\end{displaymath}

Meanwhile, in the subspace, directions in the $a_i$ blocks contribute to the volume in the same way as previously; however, each of the sets of vectors coming from one block of size $b_i$ contribute only $h^{b_i - k_i}$ times the original lengths. This gives us a volume of $h^{\beta}$, where
\begin{displaymath}
\begin{aligned}
 \beta &= (k\sum (a_i -1) + (k+1)\sum(b_i - k_i) + k + \sum (n_i -1)\\
  &= k(m-k) + 2k - c + \sum(b_i - k_i)
\end{aligned}
\end{displaymath}

Then the total volume distortion is $h^{\beta/\alpha}$.

Notice that the only part of the exponent that depends on the choice of vectors from blocks is the $\sum(b_i - k_i)$ term, which appears in both the numerator and denominator. Overall, the exponent will be maximized when this term is minimized. Thus, we want to pick our vectors so that we come as close as possible to using all of the vectors from any block from which we use more than one vector.  
\end{proof}

\subsubsection{Area distortion}\label{gerstenconj}

We can combine these results to answer a question of Gersten (\cite{gersten96}, p. 19): 

{\bf Theorem \ref{areadistthm}.}
The group $\Z^m$, $m \geq 3$, is area undistorted in $\Gamma_M$ if and only if $M$ has finite order.

See Figure~\ref{areaflow} for a flow chart for the various possible area distortion functions.

\begin{proof}
The latter condition is equivalent to saying that $M$ is diagonalizable and all eigenvalues of $M$ are roots of unity, which, by a theorem of Kronecker, is true if and only if all eigenvalues of $M$ are on the unit circle (see for example \cite{greiter78}). Then by Theorem~\ref{abcdist}, $\Z^m$ is undistorted. 

Otherwise, let us consider the possible cases. 

\begin{case}
All eigenvalues are on the unit circle.
\end{case}

Since no power of $M$ is the identity, it must be the case that some Jordan block of $M$ has ones along the superdiagonal. Then Theorem~\ref{blockdist} gives us a lower bound on area distortion.

\begin{case}
There is a block of size more than one with eigenvalue off the unit circle.
\end{case}

It will suffice in this case to show that area is distorted in the case $M = \left(\begin{array}{cc} \lambda & 1 \\ 0 & \lambda\end{array}\right)$, as this will always be a subgroup in $G_M$, giving a lower bound on the volume distortion. 

Consider the square of side length $n\lambda^{-h}$ at height $h = \frac{\log n}{\log \lambda}.$ Projecting this to height zero gives a parallelogram in $\R^2$ with area at least $n^2$; however, we can fill it in $G_M$ with five parallelograms each of area linear in $n$. Thus area is quadratically distorted.

\begin{case}
There are at least two eigenvalues off the unit circle.
\end{case}

Then by Theorem \ref{abcdist} it must be the case that the distortion is nonlinear. 

Notice that if there are at least three such eigenvalues, then two must lie on the same side of the unit circle, which means that the area distortion is quadratic, the maximum possible. 

\begin{case}
There is exactly one eigenvalue off the unit circle.
\end{case}

In this case, Theorem~\ref{oneeig} tells us the distortion is bounded below by $\frac{n^2}{W(n)}$. The distortion may in fact be higher if we have large blocks associated to unit-length eigenvectors. 

\end{proof}

In the case $m = 2$, we can also classify area distortion, though the conditions are different: here area is undistorted if and only if $\det(M) = 1$. Otherwise, $M$ has two eigenvalues, say with absolute values $\lambda$ and $\mu$, and by Theorem~\ref{abcdist}, area distortion is quadratic (maximal) if $\lambda$ and $\mu$ are both greater than one, and $n^{2 + \log_\lambda(\mu)}$ if $\lambda > 1 > \mu$ and $\lambda\mu > 1$. If $\lambda > 1$ and $\mu = 1$, then the distortion is $\frac{n^2}{W(n)}$. The examples from \S \ref{sol} and \S \ref{sull} (Sol and Nil geometry) are cases of this sort where $\det(M) = 1$.

While the cases become more complicated with higher dimensions, it should be possible to generalize Gersten's conjecture as follows.

\begin{conjecture}
Let $M$ be a square matrix of size at least $m$ with nonzero determinant, and let $k < m$. Then the group $\Z^m$ is $k$-volume undistorted in $\Gamma_M$ if and only if $M$ has finite order.
\end{conjecture}

\bibliographystyle{plain}	
\bibliography{distort}

\vspace{1cm}

\noindent
Address:\\
Department of Mathematics\\
University of Michigan\\
Ann Arbor, MI 48109\\
Email: \texttt{hbennett@umich.edu}

\end{document}